\def\tr{{\rm tr}}
\def\ax{\mbox{ax}}

\def\IH{{\Bbb H}^3}
\def\IHM{{\Bbb H}}

\def\IR{{\Bbb R}}
\def\IN{\Bbb N}
\def\IQ{\Bbb Q}

\def\ID{\Bbb D}
\def\IZ{\Bbb Z}
\def\IC{\Bbb C}
\def\IF{\Bbb F}
\def\IQ{\Bbb Q}
\def\oC{\hat{\IC}}

\documentstyle[12pt]{article}
\input{amssym.def}
\input{amssym.tex}

\title{Polynomial Trace Identities in $SL(2,{\bf C})$,  Quaternion Algebras, and Two-generator Kleinian Groups}

 \author{ T. H. Marshall \and   G. J. Martin \thanks{Research supported in
part by grants from the N.Z.
Marsden Fund.  \newline \newline
NZ Institute for Advanced Study,  Massey University,  New Zealand. \newline \newline AMS
(1991) Classification.
Primary 30C60, 30F40, 30D50, 20H10, 22E40, 53A35, 57N13, 57M60}  }
\date{}

\begin{document}

\maketitle
\newtheorem{proposition}{Proposition}[section]
\newtheorem{theorem}{Theorem}[section]
\newtheorem{lemma}{Lemma}[section]
\newtheorem{corollary}{Corollary}[section]
\newtheorem{conjecture}{Conjecture}[section]
\newtheorem{example}{Example}[section]
\newtheorem{problem}{Problem}[section]

\begin{abstract}  We study certain polynomial trace identities in the group $SL(2,\IC)$ and their application in the theory of discrete groups. We obtain canonical representations for two generator groups in \S 4 and then in \S 5 we give a new proof for Gehring and Martin's polynomial trace identities for good words,  and extend that result to a larger class which is also closed under a semigroup operation inducing polynomial composition.  This new approach is through the use of quaternion algebras over indefinites and an associated group of units. We obtain structure theorems for these quaternion algebras which appear to be of independent interest in \S 8.  Using these quaternion algebras and their units,  we  consider their relation to discrete subgroups of $SL(2,\IC)$ giving necessary and sufficient criteria for discreteness,  and another for arithmeticity \S 9.  We then show that for the groups $\IZ_p*\IZ_2$,  the complement of the closure of roots of the good word polynomials is precisely the moduli space of geometrically finite discrete and faithful representations a result we show holds in greater generality in \S12.
\end{abstract}

\newpage
\section{Introduction.}
In his work on automorphic functions \cite{FK} Fricke shows that if $\Gamma\subset SL(2,\IC)$ is a subgroup,  then the trace of any word in the generators is a polynomial with integral coefficients in the finitely many variables consisting of the traces of the generators of $\Gamma$ together with finitely many of their products.  The case $\Gamma$ has two generators,  $\Gamma=\langle A,B \rangle$ has been particulary well studied. There the trace of any word $w(A, B)\in \Gamma$ is a polynomial in the three complex variables $ x$, $y$ and $z$  where
\begin{equation}\label{xyz} x = \tr(A), \;\;\; y = \tr(B),\;\;\; z = \tr(AB).\end{equation}
As a consequence every conjugacy class of an element $w(A,B)$ in  $\Gamma$  uniquely determines a polynomial $P(x, y, z)$ with integral coefficients:    define $P$ by
\[P(x,y, z) = \tr (w(A, B)).\]
Horowitz \cite{Horo}  showed that the polynomial $P$ may not determine the conjugacy class of $w(A,B)$ uniquely, although, for a given polynomial there are only finitely many conjugacy classes represented by cyclically reduced words $w(A,B)$ giving $P$.   These results are largely based around the conjugacy invariance of trace and Fricke's simple identity
\[ \tr (AB) + \tr  (AB^{-1}) = \tr (A) \cdot \tr(B). \]
Traina \cite[Corollary 1.]{Traina} develops a family of trace identities to establish the following theorem.
\begin{theorem}[Uniqueness]
Cyclically reduced words $w_1=w_1(a,b)$ and $w_2=w_2(a,b)$ can have the same trace polynomial only if the absolute values of the exponents of the generators of $a$ in $w_2$ arise from those in $w_1$ by a permutation, and the same must be true for the exponents of $b$
\end{theorem}

\bigskip

Our initial interest lies in understanding these polynomial trace identities further and their connection with discrete groups of M\"obius transformations.  The M\"obius group acts as linear fractional transformations of the Riemann sphere;
\[PSL(2,\IC) \ni \pm  \left[ \begin{array}{cc} a & b\\ c& d \end{array}\right] \leftrightarrow \frac{az+b}{cz+d} \in \mbox{M\"ob($\oC$)} \]
and through the Poincar\'e extension we identify M\"ob($\oC$) with $Isom^+(\IH)$ the group of orientation preserving isometries of hyperbolic 3-space.  A thorough discussion of these things can be found in Beardon's book \cite{B}.  The numbers $x,y,z$ defined above at (\ref{xyz}) are not well defined in $PSL(2,\IC)$ and so we first identify new parameters.

Given two matrices $A,B\in SL(2,\IC)$ we define the parameters
\begin{equation}
\label{defineparameters}
\gamma(A,B) = \tr [A,B]-2, \hskip10pt \beta(A)=\tr^2(A)-4,\;\;  \beta(B)=\tr^2(B)-4.
\end{equation}
Here $[A,B]=ABA^{-1}B^{-1}$ is the multiplicative commutator.  These parameters depend only on the conjugacy class of $\langle A,B\rangle$ and are well defined in the projective group $PSL(2,\IC)$.  They determine the group $\langle A,B \; \big| \cdots \rangle $ uniquely up to conjugacy if $\gamma(A,B)\neq 0$, \cite{GM3}.

Note that $\gamma(A,B)$ is unchanged by Nielson moves (automorphisms of the free group of rank $2$)  on the generating pair $\{A,B\}$,  so for instance
\[ \gamma(A,B)=\gamma(B,A)=\gamma(A,A^m B^{-1} A^n) \]
and so forth.  The parameters are set up so that if $\langle A,B\rangle={\bf I}$, the trivial group, then $(\gamma,\beta,\tilde{\beta})=0\in\IC^3$.

\bigskip

In this article we are primarily interested in a special family of words,   called {\em good words},  and the family of trace polynomials they generate.   These words are defined in \S \ref{goodwords} below.  This remarkable family ${\cal W}$ has the following properties reminiscent of the Chebyshev polynomials:
\begin{enumerate} 
\item {\bf [Semigroup structure]} ${\cal W}$ forms a semigroup under the operation
\begin{equation}
w_1(a,b)*w_2(a,b) = w_1(a,w_2(a,b))
\end{equation}
\item {\bf [Polynomials and composition]} For each $w(a,b)\in {\cal W}$ there is an associated monic polynomial with integer coefficients in {\em two} complex variables $P_w(\gamma,\beta)$.  These polynomials have the property that if $w_1(a,b), w_2(a,b)\in {\cal W}$,  then
\begin{equation}
P_{w_1*w_2}(\gamma,\beta) = P_{w_1}(P_{w_2}(\gamma,\beta),\beta)
\end{equation}
That is the semigroup operation above induces polynomial composition.  
\item {\bf [Commutators and bounded roots]}  Given a representation of 
\[ \Gamma = \langle a,b \; \big| \; b^2=1 \rangle
\] into $SL(2,\IC)$,  $a\mapsto A$,  $b\mapsto B$,    set $\gamma=\gamma(A,B)$ and $\beta=\beta(A)$.  Then for $w=w(a,b)\in \Gamma$,  $w$ is a good word and we have
\begin{enumerate}
\item Commutator independence from the third complex variable $\beta(B)$, 
\[ \gamma(A,w(A,B)) = P_w(\gamma,\beta)\] 
\item Suppose $A$ is not an irrational rotation,  equivalently
\[ \beta \not\in \{-4\sin^2(r\pi) :r \in \IR \setminus  \IQ\} ,\]
 and denote the zero set of the polynomials by
\[ {\cal Z}_\beta = \{z\in \IC:P_w(z,\beta)=0 \mbox{  for some $w\in \Gamma$} \} \]
   Then $\overline{{\cal Z}} $ is compact and
\[ \IC \setminus \overline{{\cal Z}_\beta} =  {\cal R}_\beta \]
where ${\cal R}_\beta$ is nonempty,  unbounded and conformally equivalent to the punctured disk.
\item The group $\langle A,B\rangle$ is discrete and free on generators 
\[ \langle A,B\rangle\cong \langle A\rangle *\langle B\rangle\]
 if and only if $\gamma(A,B)\in \overline{{\cal R}_\beta}$.
\end{enumerate}
\end{enumerate}
Here we will prove the density of the roots of good word polynomials  in the exterior of the moduli space of discrete and faithful representations of $\Gamma$,  that is 3 (b). The hard part of 3 (c) concerns the structure of the boundary,  and the only proof we have relies on some very deep results concerning the geometry of discrete groups such as the density and the ending lamination theorems,  see \cite{NS}.  This is because ${\cal R}_\beta$ can be identified with the (moduli) space of discrete and faithful geometrically finite representations of $\IZ*\IZ_2$,  with the generator of   $\IZ=\langle A \rangle$ and $\beta(A)=\beta$.  The ``pants'' decomposition of a geometrically finite Riemann surface with fundamental group $\IZ*\IZ_2$ or $\IZ_p*\IZ_2$ shows that ${\cal R}_\beta$ is topologically a punctured disk.  These are obtained from the disk with two cone points of order two glued along its boundary to a disk with two holes (or punctures). In the case $\beta=0$,  $\overline{\cal R}$ is known in the literature as the Riley slice and the boundary $\partial {\cal R}$ is a topological circle \cite{ASWY}.  This is expected to persist for all other $\beta \in \IC$ as well.  The geometrically infinite faithful representations lie in the continuum (topological circle ?) $\partial {\cal R}$.

\bigskip

The complement $\IC\setminus \overline{{\cal R}_\beta}$ consists of nondiscrete groups apart from a countable discrete set in ${\cal R}_\beta$ of points which are the roots of  polynomials corresponding to relators in groups which are discrete but not splitting.  There are some conjectures about the structure of the polynomials and the words they come from - basically that they are associated with Dehn surgeries on two bridge knots and links and associated Hecke groups (obtained by adding an unknotting tunnel). The cusp points on the boundary arise from pinching a geodesic (arising  from a Farey word) of the Riemann surface with fundamental group $\IZ*\IZ_2$ as in \cite{KS},  giving a ray in the unbounded region ending on $\partial {\cal R}$,  while from the bounded region $\IC\setminus \partial {\cal R}$ cusps are associated with Dehn surgery limits (via Thurston's Dehn Surgery theorem) from the inside.

\medskip

Indeed it is the strong connection between these representation spaces of discrete groups,  low dimensional hyperbolic geometry and topology and the good word polynomials that motivates our consideration of them.  The geometry of commutators plays an important role in understanding the geometry and topology of discrete groups and their associated quotients, hyperbolic $3$-manifolds and $3$-orbifolds.  For instance if $A,B\in \Gamma$ where $\Gamma$ is a discrete subgroup of $SL(2,\IC)$, we put $\beta=\beta(A)$ and then suppress it writing $P_w(\gamma)$ for $P_w(\gamma,\beta)$ to find that $\{P_w(\gamma):w\in {\cal W}\}$ is a collection of traces of commutators in $\Gamma$.  Further,  if $w\in {\cal W}$, then the semigroup operation gives $w*w*\cdots*w\in {\cal W}$ and
\[ P_w(\gamma), P_{w*w}(\gamma)=P_w(P_w(\gamma)), \ldots, P_{w*w*\cdots*w}(\gamma)=P_w^{\circ n}(\gamma) \]
gives a sequence of commutator traces from the holomorphic dynamical system given by iteration of the polynomial $P_w$.  As perhaps the simplest nontrivial example,  with $w(a,b)=bab^{-1}$,  we have $P_w(\gamma,\beta)=\gamma(\gamma-\beta)$.  If $\beta=0$,  then we see $\gamma,\gamma^2,\ldots,\gamma^n\ldots$ is a sequence of commutator traces.  If $0<|\gamma|<1$, then the sequence $\{\gamma^n\}_{n\geq0}$ accumulates on $0$.  It is not a particularly difficult exercise to show this can't happen in a discrete group,  and we therefore obtain the classical Shimitzu-Leutbecher inequality.

\begin{theorem}[Shimitzu-Leutbecher inequality]\label{SL}
If $\langle A,B\rangle\subset SL(2,\IC)$ is discrete and $A$ is parabolic ($\beta=0$),  then $|\gamma|=|\tr[A,B]-2|\geq 1$. 
\end{theorem}

 J\o rgensen's inequality \cite{Jorg} follows  in the same way if $|\beta|<1$ for then $0$ is an attracting fixed point for the iterates of $P_w$ and the disk $\ID(0,1-|\beta|)$ lies in the Fatou set so $\gamma\not\in\ID(0,1-|\beta|)$ and so $|\gamma|+|\beta|\geq 1$.  We will give other examples later.

\medskip 

In order to  fully exploit these polynomials in low dimensional topology and geometry,  it is crucial to understand more about them and develop a systematic approach to uncovering the inequalities and regions of moduli space where their roots lie.  For instance,  to understand and extend the important $\frac{1}{2}\log 3$ theorem of Gabai,  Meyerhoff and Thurston \cite{GMT},  used to prove the topological rigidity of hyperbolic three manifolds \cite{Gab},  an ad hoc approach required rigorous estimates on the computation of 100+ matrix multiplications - these words were called {\em killer words} as they removed small regions of moduli space using discreteness criteria such as J\o rgensen's inequality or other criteria such as contradicting a choice of shortest geodesic.  An approach based on good words is far simpler since estimates are required for the roots of a polynomial equation with integer coefficients of lesser degree.  Such searches have been used to resolve a number of problems such as:
\begin{enumerate}
\item The unique minimal volume 3-orbifold (co-volume lattice of hyperbolic isometries) identified as the arithmetic Coxeter reflection group 3-5-3,  extended by the order two symmetry induce from the diagram, \cite{GMannals,MM}.
\item Structure of the singular set. Tables 6-10 of \cite{GMMR,GMannals} give sharp bounds for the distance between components of the singular set of a hyperbolic 3-orbifold and the distance between tetraheral, octahedral and icosahedral points in a Kleinian group.
\item  Automorphism groups of 3-manifolds and $3$-dimensional Hurwitz groups.  Sharp bounds for the order of the automorphism group of a  hyperbolic 3-manifold group in terms of the volume and analogous to the $84(g-1) $ Theorem of Hurwitz, \cite{conder}
\item Margulis constant. The Margulis constant is achieved in a two- or three- generator group,  the case of two generator groups is completely resolved, \cite{GMMargulis} and the only remaining case concerns Kleinian groups generated by three elements of order two.  
\item Geodesic length spectrum of 3-folds.  Inequalities are used to find bounds on the length of intersecting closed geodesics, or non-simple geodesics, which are within a factor of 2 of being sharp. These, together with estimates on the Margulis constant, yield good bounds for the thick and thin decompositions of hyperbolic $3$-manifolds.
\end{enumerate}

\bigskip

In this paper we uncover a group ${\cal V}$ of elements of unit norm in a quaternion algebra ${\cal Q}$ with associated indeterminates, which maps under an ``evaluation homomorphism'' $\rho:{\cal V} \to PSL(2,\IC)$
to a group which includes these good words on two
generators.  Further ${\cal V}$ naturally extends to a larger group which gives a corresponding extension of the isometry group $\rho({\cal V})$.  Roughly,  polynomials  $R,S,T$ and $W$ in the indeterminates $u$ and $v$ form a ``quaternion'' $(R,S,T,W) \in {\cal Q}$ which has norm 1 when
\begin{equation}
\label{symmdet1}
R^2-(u^2-1)S^2-(v^2-1) T^2 +(u^2-1)(v^2-1)W^2=1
\end{equation}
A special case of interest occurs when $u$ or $v$ is $\pm 1$, in which case this equation reduces to the polynomial Pell equation,
\begin{equation}
P^2(x)-(x^2-1)Q^2(x)=1,
\end{equation}

An obvious similarity between the two equations is  that the solution sets have a natural group structure,  this is what we will exploit to begin to understand the structure of good words. However,  there are significant differences.  For instance while the solutions $P(x)$, $Q(x)$ of the polynomial Pell equation must have integer coefficients, there are members of ${\cal V}$ whose polynomials have coefficients which need not even be rational (see Section \ref{irrational} for an example). We also note that (\ref{symmdet1}) has some solutions with strictly complex coefficients: a simple example is $R(u,v)=uv$, $S(u,v)=T(u,v)=1$, $W(u,v)=i$. However when we confine ourselves to solutions with rational coefficients some remarkable properties emerge; in particular it turns out all that such solutions actually have half-integer coefficients \S 8. 

\bigskip

In order to study these word polynomials more fully,  as well as justify the sorts of results we are seeking,  we need to develop a few ideas from hyperbolic geometry and in particular from the geometry of discrete groups of hyperbolic isometries of  hyperbolic $3$-space.   

\section{Background in hyperbolic geometry.}
\label{motivation}
Let ${\rm Isom}^+(\IH)$ be the group of orientation preserving isometries of $\IH$,  hyperbolic $3$-space,
\[\IH=\{x=(x_1,x_2,x_3)\in\IR^3,  x_3>0\}, \mbox{ with metric } ds=\frac{|dx|}{x_3} \]
of constant negative curvature equal to $-1$.

\medskip

We briefly review some well known facts about the group ${\rm Isom}^+(\IH)$; see e.g. \cite{B,GM3} or \cite{MM2} for more details.

\medskip

Each $f \in {\rm Isom}^+(\IH)$
is the Poincar\'e extension of a M\"{o}bius transformation of the boundary $\partial \IH$ which we identify as $\hat{\IC}=\IC\cup\{\infty\}$,  the Riemann sphere. Hence there is a natural isomorphism between ${\rm Isom}^+(\IH)$ and $PSL(2, \IC)$. 
Using the definition at (\ref{defineparameters}) we can thus define the trace and $\beta$ and $\gamma$ parameters for isometries $f, g \in Isom^+(\IH)$, simply by setting $\tr(f)=\tr(A)$, $\beta(f)=\beta(A)$ and $\gamma(f,g)=\gamma(A,B)$, where $A,B \in PSL(2, \IC)$ represent $f$ and $g$ respectively.

\medskip

Each non-identity $f \in {\rm Isom}^+(\IH)$ has either one or two fixed-points on the boundary $\hat{\IC}$. If there is just one,  then $f$ is called {\em parabolic}; if there are two, then we define the {\em axis} of $f$, $\ax(f)$ to be the hyperbolic geodesic line joining them.  Now $f$ leaves $\ax(f)$ invariant, and its action on this geodesic is a translation along by a distance $\tau=\tau(f) \geq 0$, the {\em translation length} of $f$, together with a rotation through an angle $\eta=\eta(f)$, the {\em holonomy} of $f$ around $\ax(f)$. If $\tau(f) >0$, then $\eta \in (-\pi,\pi]$, is taken anticlockwise around $\ax(f)$, as determined by the direction of the translation of $\ax(f)$ performed by $f$, and the right-hand rule; in this case $f$ is called {\em loxodromic}. If $\tau(f) =0$, that is if $f$ fixes $\ax(f)$ pointwise, then $f$ is called {\em elliptic}, in which case the distinction between clockwise and anticlockwise disappears, and we may assume that $\eta \geq 0$, that is $\eta \in (0,\pi]$. 

When $f$ is elliptic or loxodromic the parameters $\tau(f)$ and $\eta(f)$ together determine $f$ up to conjugacy. 

When $f$ is parabolic or the identity, we set $\tau(f)=\eta(f)=0$. 

\medskip

The following lemma classifies the isometries in ${\rm Isom}^+(\IH)$ up to conjugacy, and identifies, for each isometry, the conjugations which leave it unchanged.
\begin{lemma}
\label{canonicalconjugate}
A non-identity isometry $f \in {\rm Isom}^+(\IH)$ is conjugate to $z+1$ if $f$ is parabolic, and otherwise to a unique isometry of the form $f(z)=re^{i\theta}z$, where $r=e^{\tau(f)} \geq 1$, $-\pi < \theta \leq \pi$ if $r>1$, and $0 \leq \theta \leq \pi$ if $r=1$.

If $gfg^{-1}=f$, then either $g$ is the identity, $g$ and $f$ have exactly the same fixed points on $\hat{\IC}$, or, $f$ is an elliptic of order 2, and $g$ is an elliptic of order 2 which interchanges the endpoints of $\ax(f)$.

\end{lemma}

As previously remarked, both parameters $\beta(f)$ and $\gamma(f,g)$ are invariant under conjugacy. Conversely, if $\beta(f)\neq 0$,  then $\beta(f)$ determines $f$ up to congugacy,
 and if $     \gamma(f,g) \neq 0$, then $\beta(f)$, $\beta(g)$ and $\gamma(f,g)$ together determine the group $\langle f,g \rangle$ up to congugacy \cite{GM3}.  We prove this result in  Theorem \ref{canonicalgroup} below by identifying a canonical representation.
 
 \medskip 
Both the parameters $\gamma(f,g)$ and $\beta(f)$ encode geometric information.  For instance:
\begin{equation}
\label{betageom}
\beta (f)=4\sinh^2 \left(\frac{\tau +i\eta}{2}\right),
\end{equation}
and, when $f$ is elliptic or loxodromic,
\begin{equation}
\label{gammageom}
\gamma(f,g)=\frac{1}{4}\beta (f) \beta (g)\sinh^2(\Delta),
\end{equation}
where $\Delta=\Delta({\rm ax}(f),{\rm ax}(gfg^{-1}))$ represents the {\em complex distance} between $\ax(f)$ and $\ax(gfg^{-1})$ (the imaginary part of this distance, which represents the angle between the two axes, is defined modulo $\pi$, so the right hand side of (\ref{gammageom}) is well defined).
  It is an elementary fact (see e.g. \cite{B} or Theorem \ref{canonicalgroup} below) that $\gamma(f,g)=0$ if and only if $f$ and $g$ share a fixed point on the boundary $\hat{\IC}$ of $\IH$; indeed for non-parabolic $f$ and $g$, this follows immediately from (\ref{gammageom}), $\Delta$ being 0 when the axes of $f$ and $gfg^{-1}$ either meet at a point of $\hat{\IC}$ or coincide. In applications we often want to
  distinguish between these two cases. We develop an algebraic test in Section \ref{discrete}.

\section{Matrix Identities.}

We collect some matrix identities for later use.
Let

\begin{eqnarray}\nonumber
M=\left[\begin{array}{cc}
k & m \\ 0 & k^{-1}
\end{array}
\right],\;
P=\left[\begin{array}{cc}
1 & 1 \\ 0 & 1
\end{array}
\right],\; \\ \label{matrixlist}\\
Q=\left[\begin{array}{cc}
0 & i\sqrt{k} \\ i/ \sqrt{k} & 0
\end{array}
\right],\;
N=\left[\begin{array}{cc}
a & b \\ c & d
\end{array}
\right],\nonumber
\end{eqnarray}
where $ad-bc=1$. Then

\begin{equation}
\label{offdiag}
QNQ^{-1}=\left[\begin{array}{cc}
d & kc \\ b/k & a
\end{array}
\right],
\end{equation}

\begin{equation}
\label{generalconjuation}
MNM^{-1}=\left[\begin{array}{cc}
a+k^{-1}mc & -m^2c+mk(d-a)+k^2b \\ k^{-2}c & d-k^{-1}mc
\end{array}
\right],
\end{equation}
and when $m=0$

\begin{equation}
\label{generalcommutator}
[M,N]=
\left[\begin{array}{cc}
ad-k^2bc & ab(k^2-1) \\ cd(k^{-2}-1) & ad-k^{-2}bc,
\end{array}
\right]
\end{equation}

\begin{equation}
\label{paraboliccommutator}
[P,N]=
\left[\begin{array}{cc}
1+c^2+ac & 1-a^2-ac \\ c^2 & 1-ac
\end{array}
\right].
\end{equation}
In particular we have the useful trace identities, when $m=0$
\begin{equation}
\label{generaltrace}
{\rm tr}[M,N]=2-{(k-k^{-1})}^2bc
\end{equation}
\begin{equation}
\label{parabolic trace}
{\rm tr}[P,N]=2+c^2
\end{equation}

\section{Two-generator Groups.}

We now classify up to conjugacy all two-generator subgroup of ${\rm Isom}^+(\IH)$, by finding a canonical representative for each conjugacy class. Throughout we always use the principle values of square roots.
\begin{theorem}
\label{canonicalgroup}
Every group generated by two non-identity isometries in ${\rm Isom}^+(\IH)$ is conjugate to a group of the form $\langle f,g \rangle$, where $f$ and $g$ have matrix representatives $A$ and $B$ respectively in $PSL(2, \IC)$ such that either:

\medskip

\noindent{\bf Case 1.}
\begin{equation}
\label{AB}
A = \left[\begin{array}{cc}
\frac{\sqrt{\beta(f)(\beta(f)+4)}+\beta(f)}{2\sqrt{\beta(f)}}
 & 0 \\ 0 & \frac{\sqrt{\beta(f)(\beta(f)+4)}-\beta(f)}{2\sqrt{\beta(f)}}
\end{array}
\right],\;\;  B=\left[\begin{array}{cc}
a & b \\ c & d
\end{array}
\right]  
\end{equation}
\begin{eqnarray}
a&=&
\frac{1}{2}\left(\sqrt{\beta(g)+4}+\sqrt{\frac{
4\gamma(f,g)+\beta(f)\beta(g)}{\beta(f)}}
\right) \nonumber \\
\label{ad}
d &=&\frac{1}{2}\left(\sqrt{\beta(g)+4}-\sqrt{\frac{
4\gamma(f,g)+\beta(f)\beta(g)}{\beta(f)}}
\right),
\end{eqnarray}
and
\begin{equation}
\label{gammanonzerooffdiagonals}
c=-b= \sqrt{\frac{\gamma(f,g)}{\beta(f)}}
\end{equation}
when $\gamma(f,g) \neq 0$,  and either
\begin{equation}
\label{offdiagonals}
b=0,c=1\;\;\mbox{or}\;\; b=1, c=0\;\;\mbox{or}\;\; b=c=0
\end{equation}
when $\gamma(f,g) =0$;
or \\

\medskip

\noindent{\bf Case 2.}
\begin{equation}
\label{par1}
A=\left[\begin{array}{cc}
1 & 1\\ 0 & 1
\end{array}
\right],\;\;\;B=\left[\begin{array}{cc}
0 & -1/\sqrt{\gamma(f,g)}\\ \sqrt{\gamma(f,g)} &  \sqrt{\beta(g)+4}
\end{array}
\right],
\end{equation}
or \\

\medskip

\noindent{\bf Case 3.}
\begin{eqnarray}
\nonumber
A& = & \left[\begin{array}{cc}
1 & 1\\ 0 & 1
\end{array}
\right], \\ B &=&
\left[ \begin{array}{cc} \frac{1}{2}[\sqrt{\beta(g)+4}+\sqrt{\beta(g)}] & \ell \\ 0 & \frac{1}{2}[\sqrt{\beta(g)+4}-\sqrt{\beta(g)}]  \end{array} \right] \label{par2}
\end{eqnarray}
where $\ell=0$ when $\beta(g) \neq 0$, and can take any complex value when $\beta(g)= 0$.

\medskip

The three cases are respectively the cases $\beta(f) \neq 0$ ($f$ non-parabolic), $\beta(f)=0$ and $\gamma(f,g) \neq 0$ and $\beta(f)=\gamma(f,g) =0$ ($f$ parabolic).
\end{theorem}

\bigskip

{\bf Proof.} We set $\beta=\beta(f)$, $\gamma=\gamma (f,g)$. Suppose first that $f$ is loxodromic or elliptic ($\beta \neq 0$). By Lemma \ref{canonicalconjugate} we can conjugate $f$ so that its matrix representative $A$ is diagonal.

By (\ref{generaltrace}) we have that $\gamma = -bc\beta$, whence
\begin{equation}
\label{bc}
bc=\frac{-\gamma}{\beta}, \;\;\; ad=1-\frac{\gamma}{\beta}
\end{equation}
We have ${(a+d)}^2-4=\beta(g)$, and since $B$ is determined only up to sign, we may thus assume that $a+d=\sqrt{\beta(g)+4}$. Together with (\ref{bc}), this gives that $a$ and $d$ are either as given by (\ref{ad}),
or are obtained from these by interchanging the values of $a$ and $d$.
Using (\ref{offdiag}), we may then conjugate $A$ and $B$ if necessary, to interchange $a$ and $d$, so that (\ref{ad}) holds, and $A$ is still diagonal.

Let $r$ and $s$ be the diagonal entries of $A$.
We have ${(r+s)}^2-4=\beta$, and since $A$ is determined only up to sign, we may assume that $r+s=\sqrt{\beta(\beta+4)}/\sqrt{\beta}$. Together with the condition $rs=1$, this gives that either $A$ or $A^{-1}$ takes the form given by (\ref{AB}). Since $\langle A^{-1}, B \rangle$ and $\langle A, B \rangle$ are the same group, we may assume that
$A$ satisfies (\ref{AB}). 

\medskip

Finally, we apply a conjugacy of the type (\ref{generalconjuation}) (with $m=0$) to $A$ and $B$ to adjust the values of $b$ and $c$, leaving $A$ unchanged. If $b,c \neq 0$ (i.e. when $\gamma\neq 0$), we can use such a conjugacy to give $b$ are $c$ any values subject to (\ref{bc}); in particular we can make (\ref{gammanonzerooffdiagonals}) hold. If exactly one of the values of $b$ and $c$ is nonzero, then we conjugate to make it 1. The only other possibility is $b=c=0$, so the options given in (\ref{gammanonzerooffdiagonals}) and (\ref{offdiagonals}) are exhaustive.

Now we suppose that $f$ is parabolic. Using Lemma \ref{canonicalconjugate} we conjugate so that $f(z)=z+1$, so that its matrix representative $A=\left[\begin{array}{cc}
1 & 1\\ 0 & 1
\end{array}
\right]$.
By (\ref{parabolic trace}) $c^2=\gamma$, and we may assume, since $B$ is determined only up to sign, that $c=\sqrt{\gamma}$.

Now we have two subcases, determined by whether or not $\gamma=0$.
If $\gamma \neq 0$, then $c \neq 0$, and we apply the conjugation (\ref{generalconjuation}) with $k=1$ and $m=-a/c$ to $A$ and $B$, leaving $A$ the same, and changing $a$ to 0, whence
\begin{equation}
\label{conjevenpower}
B=\left[\begin{array}{cc}
0 & -1/\sqrt{\gamma}\\ \sqrt{\gamma} & \pm \sqrt{\beta(g)+4}
\end{array}
\right].
\end{equation}
Let $B_{(+)}$ and $B_{(-)}$ be the matrices obtained by taking the $+$ and $-$ signs respectively in (\ref{conjevenpower}). We show that $\langle A, B_{(+)} \rangle$ and $\langle A, B_{(-)} \rangle$ are conjugate groups in $PSL(2, \IC)$. A conjugation of the form (\ref{generalconjuation}) with $k=i$ and $m=0$ takes
$A$ to $A^{-1}$ and $B_{(-)}$ to $-B_{(+)}$. Thus $\langle A, B_{(-)} \rangle$ is conjugate to
$\langle A^{-1}, -B_{(+)} \rangle=\langle A, B_{(+)} \rangle$ in $PSL(2, \IC)$ as required. Thus, without loss of generality, we take the $+$ sign in (\ref{conjevenpower}).

If $\gamma=0$, then $c=0$ so $ad=1$. As in previous cases, we may assume that $a+d=\sqrt{\beta(g)+4}$, so that $B$ must take the form (\ref{par2}), up to an interchange of the diagonal entries. If $\beta(g)=0$, then these entries are the same, and we are done. Otherwise $a \neq d$, and we can apply a further conjugation  of the form (\ref{generalconjuation}), with $k=1$, so as to get both $\ell=0$ and to leave $A$ unchanged. Now $B$ is diagonal, and interchanging $a$ and $d$ replaces $B$ by $B^{-1}$. Since this operation leaves the group $\langle A,B \rangle$ unchanged, we may assume that $B$ is given by (\ref{par2}). \hfill
\hfill $\Box$

\bigskip

\noindent {\bf Remarks.}
We can characterize geometrically the four ways of assigning values to $b$ and $c$ given by (\ref{gammanonzerooffdiagonals}) and (\ref{offdiagonals}). In (\ref{gammanonzerooffdiagonals}) $\gamma \neq 0$, so $b,c \neq 0$ and $f$ and $g$ have no common fixed points in $\hat{\IC}$. In this case, as remarked in the proof, $b$ and $c$ can be made to take any values whose product is $-\gamma/\beta$, and the exact choice is rather arbitrary.
However the normalization that we have chosen is quite natural from a geometric viewpoint; it makes the fixed points of $g$ mutually reciprocal, and (consequently), when $g$ is non-parabolic, the common perpendicular of $\ax(f)$ and $\ax(g)$ is the geodesic with endpoints $\pm 1$. When $g$ is parabolic the fixed point is $z=1$. See \cite{MM} for more details.

 The first two cases of (\ref{offdiagonals}), when $\{b,c\}=\{0,1\}$
 occur when $f$ and $g$ have a single common fixed point in $\hat{\IC}$. If $f$
is loxodromic, then this point is repulsive when $b=0$, $c=1$, and attractive when
 $b=1$, $c=0$.
  If $f$ is elliptic, then $f$ rotates $\IH$ anticlockwise (resp. clockwise) around $\ax(f)$ oriented away from the shared fixed point when $b=0$, $c=1$ (resp. $b=1$, $c=0$). (When the elliptic $f$ is order two these two cases are conjugate.) Finally $b=c=0$ when $f$ and $g$ are both elliptic or loxodromic, and have the same axis.

Note that, when $\gamma=0$, although $\beta(f)$, $\beta(g)$ and $\gamma$ do not determine the conjugacy class of $\langle f,g \rangle$,
when $\beta(f) \neq 0$ (and symmetrically when $\beta(g) \neq 0$) then there there are only three possibilities. Only when $\beta(f)=\beta(g)=\gamma=0$ (Case 3 of the theorem with $\beta(g)=0$) do the same parameters give an infinite family of non-conjugate groups.

\section{Good words} \label{goodwords}
A {\em good word} on the letters $a$ and $b$ is a word of the form
\[
w(a,b)=b^{s_1} a^{r_1} b^{s_2} a^{r_2} \ldots  b^{s_{m-1}}  a^{r_{m-1}}  b^{s_m}
\]
where $s_1=\pm 1$,
$s_j={(-1)}^{j+1}s_1$, and the $r_j$ take integer values.  

\medskip

Thus the powers of $b$ in a good word alternate in sign. By setting $r_1=0$  (resp. $r_{m-1}=0$), we obtain a good word which begins (resp. ends) with a power of $a$.

\medskip

A good word is {\em even} if
$r_1+r_2+\ldots r_{m-1}$ is even, {\em odd} otherwise, {\em balanced} if
$m$ is even, {\em unbalanced} otherwise and {\em regular} if $s_1=1$, {\em irregular} otherwise.
If $r_j=0$ for any $1<j<m-1$, then $w(a,b)$ collapses into a shorter good word (which has the same balance, parity and regularity as the original word), so we may assume that these interior powers are non-zero.

\medskip

The following easy observation is quite useful.  We leave the proof to the reader.
\begin{theorem}
\label{groupindex}
The regular balanced words in $\langle a,b\rangle$,  say $\Gamma_{reg}$,  comprise a subgroup of the free group on $a$ and $b$, of which the regular balanced even words form an index-two subgroup.
\end{theorem}

\begin{lemma}
\label{evengen}
The group of regular balanced even words,  say $\Gamma_{reg}^{even}$,  on $a$ and $b$ is generated by $a^2$, $ba^2b^{-1}$ and $[b,a]=bab^{-1}a^{-1}$.
\end{lemma}
{\bf Proof.} Let
$w=b a^{r_1} b^{-1} a^{r_2} \ldots  b  a^{r_{2m-1}}  b^{-1}$ be regular, balanced and even. We use
induction on $m$. We have

\begin{equation}
ba^ib^{-1}a^j={(ba^2b^{-1})}^{i/2}{(a^2)}^{j/2}
\end{equation}
when $i$ and $j$ are both even, and
\begin{equation}
ba^ib^{-1}a^j={(ba^2b^{-1})}^{(i-1)/2}[b,a]{(a^2)}^{(j+1)/2}
\end{equation}
when $i$ and $j$ are both odd.
This deals with the case $m=1$ and the induction step when $r_1$ and $r_2$ have the same parity.

 If $r_1$ and $r_2$ have the opposite parity, then $m>1$, and we use the same identities together with
\[
b a^{r_1} b^{-1} a^{r_2} b a^{r_3} b^{-1} a^{r_4}=(b a^{r_1} b^{-1} a^{r_2-1}){[b,a]}^{-1}
(b a^{r_3+1} b^{-1} a^{r_4})
\]
This completes the proof. \hfill $\Box$

\medskip
The next corollary is also immediate.

\begin{corollary}
Suppose that $a$ has order three,  $a^3=1$.  Then the group of regular balanced even words on $a$ and $b$ is a two generator group generated by $a$ and $bab^{-1}$.
\end{corollary}

We recall  here the well known identity
\begin{eqnarray*}
[b,a]& = &  (ba)^2   \; (a^{-1} b^{-1} a)^2 \; (a^{-1})^2 
\end{eqnarray*}
This tells us that the regular balanced words in $\Gamma=\langle a,b \rangle$ lie in the group $\Gamma^{(2)}$ generated by squares of elements.

\begin{corollary}
\label{evengen}
The group of regular balanced even words on $a$ and $b$ lies in the group generated by the four squares
\[  \Gamma_{reg}^{even} < \langle a^2, (bab^{-1})^2,  (ba)^2, (aba^{-1} )^2 \rangle  \]
\end{corollary}
The remark following the next result,  Theorem \ref{matrixrep},  shows that the index between these two groups is infinite.  In fact for any representation in $SL(2,\IC)$ the trace fields are $\IQ(\tr \Gamma_{reg}^{even} )=\IQ(\beta,\gamma)$ and $\IQ(\beta,\gamma, \beta(ba)) $,  using our earlier notation.

\bigskip

We can now state our first main theorem.
\begin{theorem}
\label{matrixrep}
Let $w=w(a,b)=b a^{r_1} b^{-1} a^{r_2} \ldots  b  a^{r_{2m-1}}  b^{-1}$ be a regular balanced even word, then there are polynomials $r_w,s_w,t_w,w_w$, such that
\begin{equation}
\label{integercoeffs}
2r_w,2s_w,2t_w,2w_w, r_w-s_w ,t_w-w_w \in \IZ[x,z],
\end{equation}
$r_w(0,0)=1$ and
\[ g_w(x,z):=\frac{s_w(x,z)-zw_w(x,z)}{x} \]
 is also a polynomial, and if
$f,g \in {\rm Isom}^+(\IH)$ are not the identity, have parameters $\beta=\beta(f)$, $\beta'=\beta(g)$ and $\gamma=\gamma(f,g)$, and if $A$ and $B$ are the matrices from Theorem \ref{canonicalgroup} which represent (up to conjugacy) $f$ and $g$ respectively, then:
for $f$ non-parabolic ($\beta \neq 0$)
\begin{equation}
\label{matrixform1}
w(A,B)=\left[\begin{array}{cc}
r_w(\beta, \gamma)+\frac{s_w(\beta, \gamma)Q}{\beta} & ab[\beta t_w(\beta, \gamma)+w_w(\beta, \gamma)Q]  \\ cd[\beta t_w(\beta, \gamma)-w_w(\beta, \gamma)Q] & r_w(\beta, \gamma)-\frac{s_w(\beta, \gamma)Q}{\beta}
\end{array}
\right]
\end{equation}
where $Q=\sqrt{\beta(\beta+4)}$ and $a$, $b$, $c$, $d$ are as in (\ref{ad}) and (\ref{gammanonzerooffdiagonals});
and for $f$ parabolic ($\beta = 0$), with $\gamma \neq 0$
\begin{equation}
\label{evzero1}
w(A,B)=\left[\begin{array}{cc}
r_w(0, \gamma)+\gamma t_w(0, \gamma) & 4g_w(0, \gamma)+2w_w(0,\gamma)  \\ 2\gamma w_w(0,\gamma) & r_w(0, \gamma)-\gamma t_w(0, \gamma)
\end{array}
\right],
\end{equation}
and for $f$ parabolic, with $\gamma = 0$
\begin{equation}
\label{evzerozero1}
w(A,B)=\left[\begin{array}{cc}
1 & 4g_w(0,0)-\left(\beta'+\sqrt{\beta'}\sqrt{\beta'+4}\right)w_w(0,0)
 \\ 0 & 1
\end{array}
\right]
\end{equation}

In particular, the trace ${\rm tr} (w(f,g))=2r_w(\beta, \gamma)\in \IZ[x,z]$.
\end{theorem}

\noindent {\bf Remark.}
A key feature here is that the polynomials $r_w$, $s_w$, $t_w$ and $w_w$, and in case (\ref{evzero1}) the whole matrix, are independent of  $\beta'=\beta(g)$. In particular this is true of traces of the matrix representations above.

\bigskip

We prove the above result in Section \ref{quat}. We can also use it to find $w(A,B)$ when $w$ is unbalanced; we do this next for non-parabolic $f$.


\begin{corollary}
\label{unbalancedword}
If $w=w(a,b)=b a^{r_1} b^{-1} a^{r_2} \ldots  b^{-1}  a^{r_{2m-2}}  b$ is a regular unbalanced even word, then there are polynomials $r_w,s_w,t_w,w_w$ with half-integer coefficients
such that for $f$ non-parabolic
\begin{equation}
\label{matrixformub}
w(A,B)=\left[\begin{array}{cc}
a(r_w(\beta, \gamma)+s_w(\beta, \gamma)Q) & b(t_w(\beta, \gamma)+w_w(\beta, \gamma)Q)  \\ c(t_w(\beta, \gamma)-w_w(\beta, \gamma)Q) & d(r_w(\beta, \gamma)-s_w(\beta, \gamma)Q)
\end{array}
\right],
\end{equation}
where $a$, $b$, $c$, $d$ are as in (\ref{ad}) and (\ref{offdiagonals}).
\end{corollary}
{\bf Proof.} First note that $\tilde{w}=w*b^{-1}$ is balanced, thus  
\[
w(A,B)=\tilde{w}(A,B)B=\left[\begin{array}{cc}
r_{\tilde{w}}+\frac{s_{\tilde{w}}Q}{\beta} & ab(\beta t_{\tilde{w}}+w_{\tilde{w}}Q)  \\ cd(\beta t_{\tilde{w}}-w_{\tilde{w}}Q) & r_{\tilde{w}}-\frac{s_{\tilde{w}}Q}{\beta}
\end{array}
\right]\left[\begin{array}{cc}
a & b \\ c & d
\end{array}
\right] = \] 
\[
\left[\begin{array}{cc}
a[r_{\tilde{w}}+\frac{s_{\tilde{w}}Q}{\beta}-\frac{\gamma}{\beta}(\beta t_{\tilde{w}}+w_{\tilde{w}}Q)[& b[r_{\tilde{w}}+\frac{s_{\tilde{w}}Q}{\beta}+\left(1-\frac{\gamma}{\beta}\right)(\beta t_{\tilde{w}}+w_{\tilde{w}}Q)]
\\ c[r_{\tilde{w}}-\frac{s_{\tilde{w}}Q}{\beta}+\left(1-\frac{\gamma}{\beta}\right)(\beta t_{\tilde{w}}-w_{\tilde{w}}Q)]
& d[r_{\tilde{w}}-\frac{s_{\tilde{w}}Q}{\beta}-\frac{\gamma}{\beta}(\beta t_{\tilde{w}}-w_{\tilde{w}}Q)]
\end{array}
\right] 
\]
\[
= \left[\begin{array}{cc}
a\left((r_{\tilde{w}}-\gamma t_{\tilde{w}})+g_{\tilde{w}}Q
\right) & b[r_{\tilde{w}}+(\beta-\gamma) t_{\tilde{w}}+(g_{\tilde{w}}+w_{\tilde{w}})Q]
\\ c[r_{\tilde{w}}+(\beta-\gamma) t_{\tilde{w}}-(g_{\tilde{w}}+w_{\tilde{w}})Q]
& d\left((r_{\tilde{w}}-\gamma t_{\tilde{w}})-g_{\tilde{w}}Q
\right) \end{array}
\right]
\]
using (\ref{bc}). We have now found the polynomials for $w$ explicitly in terms of those for $\tilde{w}$: 
\begin{eqnarray*} r_w&=&r_{\tilde{w}}-\gamma t_{\tilde{w}}\\s_w&=&g_{\tilde{w}}\\t_w&=&r_{\tilde{w}}+(\beta-\gamma) t_{\tilde{w}}\\
w_w&=&g_{\tilde{w}}+w_{\tilde{w}}
\end{eqnarray*} 
\hfill $\Box$

\bigskip

We now extend the definition of the polynomials $r_w$, $s_w$, $t_w$ and $w_w$ to arbitrary good words.

\bigskip

\noindent {\bf Definition.} Let $w$ be a good word. If $w$ is regular and even, then then the polynomials $r_w$, $s_w$, $t_w$ and $w_w$ are as defined in Theorem \ref{matrixrep} for balanced $w$, and Corollary \ref{unbalancedword} for unbalanced $w$. If $w$ is regular and odd, then $v:=w.a$ is regular and even and we define $r_w=r_v$, $s_w=s_v$, $t_w=t_v$ and $w_w=w_v$. Finally, if $w=w(a,b)$ is irregular, then $w'=w(a,b^{-1})$ is regular, and we define $r_w=r_{w'}$ etc.
Note that ${(w.a)}'=w'.a$, so that if $w$ is irregular {\em and} odd, we have $r_w=r_{{(w.a)}'}=r_{w'.a}$ and so forth.

\bigskip

If $w=w(a,b)$ is an irregular even word,
then $w'$ is regular and even.
Since $\beta(g^{-1})=\beta(g)$ and $\gamma(f,g^{-1})=\gamma(f,g)$, Theorem \ref{matrixrep} and Corollary \ref{unbalancedword} show that $w(f,g)=w'(f,g^{-1})$ is conjugate to $w'(f, g)$; in particular these have the same trace $2r_{w}(\beta, \gamma)=2r_{w'}(\beta, \gamma)$.

\medskip

Since for any good word $w=w(a,b)$, the commutator $[a,w]$ is balanced and even, it
follows easily that $\gamma(f, w(f,g))=p_w(\beta, \gamma)$, for some polynomial $p_w \in \IZ[x,z]$ as first observed in {\cite{GM3}}. The next result expresses these polynomials in terms of $t_w$ and $w_w$.
\begin{theorem}
\label{tracepolynomials}
Let $w=w(a,b)$ be a good word, $\beta=\beta(f)$, $\gamma=\gamma(f,g)$ and $h=w(f,g)$, then there is a polynomial $p_w \in \IZ[x,z]$, such that
\begin{equation}
\gamma(f,h)=p_w(\beta, \gamma)
\end{equation}
If $w$ is balanced, then
\begin{equation}
\label{evenk}
p_w(x,z)=z(x-z)(xt_w^2(x,z)-(x+4)w_w^2(x,z)),
\end{equation}
If $w$ is unbalanced, then
\begin{equation}
\label{oddk}
p_w(x,z)=z
(t_w^2(x,z)-x(x+4)w_w^2(x,z)).
\end{equation}
\end{theorem}
{\bf Proof for non-parabolic $f$.}
We may assume that $f$ and $g$ have matrix representatives $A$ and $B$ respectively as given by (\ref{AB})-(\ref{offdiagonals}). Suppose first that $w$ is regular and even, and set $W=w(A,B)$. The required identities follow from (\ref{generaltrace}) and (\ref{bc}), together with
(\ref{matrixform1}) when $W$ is balanced and (\ref{matrixformub}) when $W$ is unbalanced.
For balanced $w$, $p_w \in \IZ[x,z]$ follows from (\ref{integercoeffs}). For unbalanced $w$, using the results and notation of the proof of Corollary \ref{unbalancedword}, we have
\begin{eqnarray*}
t_w-xw_w &=& r_{\tilde{w}}+(x+z)t_{\tilde{w}}-x(g_{\tilde{w}}+w_{\tilde{w}})\\
&=& r_{\tilde{w}}+(x+z)t_{\tilde{w}}-(s_{\tilde{w}}-zw_{\tilde{w}})-xw_{\tilde{w}}\\
&=& (r_{\tilde{w}}-s_{\tilde{w}})+(x+z)(t_{\tilde{w}}-w_{\tilde{w}})\in \IZ[x,z],
\end{eqnarray*}
again by (\ref{integercoeffs}). Consequently also $t_w+xw_w \in \IZ[x,z]$, and the same then follows for $p_w$, as given by (\ref{oddk}).

If $w$ is odd, the result follows from the previous case, together with the identity $[A, W]=[A, WA]$. Finally suppose that $w$ is irregular, then the commutator $[a,w]$ is irregular and even, and
$[a,w]'=[a,w']$ is regular. Since, as previously noted, $\tr ([a,w])=\tr ([a,w]')=\tr ([a,w'])$, we are reduced to the case where $w$ is regular.
\hfill $\Box$

\bigskip

It is not difficult to prove the above directly when $f$ is parabolic. However we will instead use a limiting argument in Section \ref{quat} (Corollary \ref{limitingcases}).

\subsection{Examples of word polynomials.} \label{examples}

Before going too much further it is worthwhile giving a few examples of polynomials.  These appear in Table 1 below. 

 \medskip

{\bf Table 1. Some examples of word polynomials. }\\
{\small 
\medskip
\noindent\begin{tabular}{|c|c|c|c|c|}
\hline
 Polynomial & word    \\
 \hline
  \hline
$\gamma  (\gamma -\beta )$ & $ bab $  \\   
\hline
$(\beta +4) (\gamma-\beta ) \gamma $ & $ba^2b$ \\
\hline
$ (\beta -\gamma +1)^2 \gamma $ & $ babab$ \\
\hline
$\gamma (1-2 \beta +\gamma ^2-(\beta -2) \gamma )$ & $ baba^{-1} b$ \\
\hline 
$\gamma   (1+\beta (\beta +1) (\beta +4)-(\beta +4) (2 \beta +1) \gamma+(\beta +4) \gamma^2)$ & $ baba^2b$\\
\hline\
$\left(\beta ^2-(\gamma -4) \beta -4 \gamma +1\right)^2 \gamma $ & $ ba^2ba^2b$\\
\hline
$  \gamma(\gamma-\beta) (\beta -\gamma +2)^2 $ & $bababab$\\
\hline
$ \gamma (\beta ^2+\gamma^3-2 \beta  \gamma^2+(\beta -1) \beta  \gamma )$ & $bababa^{-1}b$\\
\hline
$ \gamma(\beta +4)  (\beta ^2+\gamma^3-2 \beta  \gamma^2+(\beta -1) \beta \gamma ) $ & $baba^2ba^{-1}b$\\
\hline
$\gamma ^3(\gamma -\beta ) (\beta +4)  (\beta  (\gamma ^2-3\gamma-4)-\beta ^2( \gamma+1) +4 \gamma ^2+4 \gamma +1)$ & $ba^{-2}bababa^{-2}bab$\\
\hline
\end{tabular}
\medskip
}
The last polynomial here is quite long,  but it has the remarkable property that it has $\gamma=0$ as a super-attracting fixed point.  We will need this fact later.

\section{Applications}
\label{discrete}
If $\langle f, g \rangle$ is discrete and non-elementary,  which in this setting the latter means that $\langle f, g \rangle$ is not virtually abelian,  then we have J\o rgensen's inequality,  \cite{Jorg}.  We sketched a proof for this earlier at Theorem \ref{SL}.
\begin{theorem} Let $\langle f,g\rangle$ be a discrete nonelementary subgroup of $SL(2,\IC)$.  Then
\begin{equation}
\label{Jorgensen}
|\beta(f)|+|\gamma (f,g)| \geq 1.
\end{equation}
\end{theorem}
 If $w$ is a good word and if $h=w(f,g)$, then we have $\gamma (f,h)=p_w( \gamma,\beta)$, where $p_w$ is the polynomial of Theorem \ref{tracepolynomials} and so we deduce from J\o rgensen's inequality applied to the group $\langle f,h\rangle$,  that
\[ |\beta(f)|+|\gamma (f,h)| = |\beta|+|p_w(\beta, \gamma)| \geq 1 \]
unless $\gamma (f,h)  = 0$, in which case $\langle f, h \rangle$ is elementary.  

\medskip

We would like to understand how this exception happens more generally.

\medskip

If $|\beta(f)|+|\gamma (f,h)| < 1$ and $\gamma (f,h) \neq 0$, then $\langle f, g \rangle$ is non-elementary since $f$ is either loxodromic or elliptic of order at least $7$.  Then J\o rgensen's inequality implies that $\langle f, g \rangle$ cannot be discrete. On the other hand, this group may be elementary if $\gamma (f,h)=0$.  However in this case we know that the fixed point sets of $f$ and of $h$ on $\hat{\IC}$ intersect, and they may coincide; for non-parabolic $f$ and $h$, coincidence means that they have the same axis. We can use polynomials to determine when this happens.

\begin{proposition}
Let $\langle f, g \rangle$ be discrete and non-elementary,   and $h=w(f,g)$ for a good word $w$.  Set $\beta=\beta(f) \neq -4$, and $\gamma=\gamma(f,g) \neq \beta, 0$.  Then the fixed point sets of $f$ and $h$  on $\hat{\IC}$ coincide if and only if $t_w(\beta, \gamma)=w_w(\beta, \gamma)=0$
\end{proposition}

\noindent {\bf Proof.} We suppose first that $w$ is regular and even. We may suppose that $f$ and $g$ are normalized so that their matrix representatives are as given in Theorem \ref{canonicalgroup} (specifically, by (\ref{AB})-(\ref{offdiagonals}) for $\beta \neq 0$, and by (\ref{par2}) for $\beta = 0$), and consequently that the matrix representative for $h$ is given by Theorem \ref{matrixrep}.

Suppose first that $\beta \neq 0$, so that $f$ is loxodromic or elliptic with fixed points $0$ and $\infty$ on the boundary. Since $abcd=\gamma (\gamma-\beta)/ \beta^2 \neq 0$, (\ref{matrixform1}) and (\ref{matrixformub}) show that $h$ shares these fixed points if and only if
\[ \beta t_w(\beta, \gamma) \pm w_w(\beta, \gamma)Q=0,\] when $w$ is balanced, and \[ t_w(\beta, \gamma) \pm w_w(\beta, \gamma)Q=0,\] when it is not. Since also $\beta \neq 0$ and $Q \neq 0$, this is equivalent to $t_w(\beta, \gamma)=w_w(\beta, \gamma)=0$ in both cases.

Now suppose first that $\beta=0$, so that $f(z)=z+1$ with fixed point $\infty$. If $\infty$ is also the only
fixed point of $h$, then $h$ must also be parabolic, and by (\ref{evzero1}), we must have $r_w(0, \gamma)=1$ and
$w_w(0, \gamma)=0$, in order to get the right trace and the right fixed point respectively. The determinant condition then gives $t_w(0, \gamma)=0$. The converse is clear.

If $w$ is regular and odd then the fixed point sets of $f$ and $h$ on $\hat{\IC}$ coincide if and only if the same is true of the fixed point sets of $f$ and $hf$. Since $w . a$ is even this occurs exactly when $t_{w . a}=w_{w.a}=0$, by the previous case, and since
$t_w=t_{w .a}$ and $w_w=w_{w . a}$, we are done.
If $w$ is irregular then $h=w'(f,g^{-1})$, so that, since $\beta(g^{-1})=\beta(g)$, $\gamma(f,g^{-1})=\gamma(f,g)$ and $w'$ is regular, the previous cases give that the fixed point sets of $f$ and $h$  on $\hat{\IC}$ coincide if and only if $t_{w'}(\beta, \gamma)=w_{w'}(\beta, \gamma)=0$, and we are done, since by definition, $t_w=t_{w'}$ and $w_w=w_{w'}$.
\hfill $\Box$

\bigskip

The above, together with (\ref{evenk}) and (\ref{oddk}), gives the following corollary which is useful.  It allows us to obtain a contradiction to discreteness unless we can identify a multiple root,  which is easy to do computationally, using the vanishing of a discriminant.

\begin{corollary}
Let $f$, $g$ and $h$ be as above and let $\beta \neq -4$, $\beta \neq \gamma$ and $\gamma \neq 0$.  If  $f$ and $h$ have the same fixed points in $\oC$, then $z=\gamma$ is a multiple root of $p_w(\beta,z)$.
\end{corollary}

\noindent {\bf Remark.} The converse of this result is false. For example, when $f$
is parabolic, $f$ and $h$ cannot have the same fixed points on the boundary unless $h$ is also parabolic, but
inspection of (\ref{evenk}) and (\ref{oddk}) shows that $p_w(0,z)$ has no simple roots apart from $z=0$.

\section{Quaternion Algebras}
\label{quat}
In order to prove Theorem \ref{matrixrep}, we switch from matrices into the slightly more abstract setting
of quaternion algebras. The objects we will be dealing with here are essentially the same as matrices of the form (\ref{matrixform1}), but with numbers replaced by indeterminates (we will let $x$ and $z$ correspond to $\beta$ and $\gamma$ respectively), and with square roots defined abstractly. We first recall some basic facts; see \cite{MR} or \cite{V} for more details.

\bigskip

For each field ${\IF}$ of characteristic $\neq 2$ and non-zero $a,b \in {\IF}$, the quaternion algebra
\begin{equation}
\label{quatalg}
{\cal A}=\left( \frac{ a, b}{\IF} \right),
\end{equation}
is defined to be the associative algebra over the field ${\IF}$ with multiplicative identity 1 and basis $\{1,{\bf i},{\bf j},{\bf k}\}$, with multiplication determined by ${{\bf i}}^2=a$, ${{\bf j}}^2=b$, ${\bf i}{\bf j}=-{\bf j}{\bf i}={\bf k}$, whence also
${{\bf k}}^2=-ab$, ${\bf j}{\bf k}=-b{\bf i}$ and ${\bf k}{\bf i}=-a{\bf j}$. The generic member of ${\cal A}$ is thus
$x+y{\bf i} +z{\bf j}+w{\bf k}$, where $x,y,z,w \in \IF$; we may abbreviate this to $(x,y,z,w)$. 

\medskip

${\cal A}$ is also isomorphic to the algebra of matrices of the form
\begin{equation}
\label{matmap}
\left[\begin{array}{cc}
x+y\xi_1 & (z+w\xi_1)\xi_2 \\ (z-w\xi_1)\xi_2 & x-y\xi_1
\end{array}
\right]
\end{equation}
over the extension field $\IF(\xi_1, \xi_2 )$, where $\xi_1^2=a$ and $\xi_2^2=b$. 

\medskip

If $n=(x,y,z,w) \in {\cal A}$, then the {\em conjugate} of $n$ is given by 
\[ \overline{n}=(x,-y,-z,-w),\] and the {\em norm} of $n$ by 
\[ N(n)=n \overline{n}=x^2-ay^2-bz^2+abw^2=x^2-ay^2-b(z^2-aw^2) \in \IF.\]
Note that the norm becomes the determinant under the mapping (\ref{matmap}).

\medskip

For reference, multiplication in the quaternion algebra (\ref{quatalg}) is given explicitly by
\begin{eqnarray}
\lefteqn{(x,y,z,w)(x',y',z',w')} \nonumber
\\ &=&(xx'+ayy'+bzz'-abww',\; xy'+yx'+b(wz'-zw'), \nonumber \\
\label{quatmult}
&& xz'+zx'+a(yw'-wy'),\; xw'+wx'+yz'-zy').
\end{eqnarray}

If the field $\IF$ is the field of fractions of an integrally closed integral domain $R$, we define an {\em $R$-lattice} in ${\cal A}$ to be a finitely generated $R$-module $L$ in ${\cal A}$; $L$ is an {\em ideal}
if $\IF L={\cal A}$. An element $\alpha \in {\cal A}$ is an integer (over $R$) if
$R[\alpha]$ is an $R$-lattice.
An {\em order} in ${\cal A}$ is an ideal which is also ring with 1. 
 By contrast with the commutative case, the set of all integers in ${\cal A}$ is not generally a ring.
 
 We will be particularly concerned with elements of norm 1. Note that these are units, and in any order which is closed under conjugation, they form a multiplicative group.

\bigskip

We consider the following quaternion algebra over the field of rational functions in two indeterminates.

\[
{\cal Q}_0:=\left( \frac{(x+4)/x, z(z-x)}{\IR(x,z)} \right)
\]

\begin{theorem}
\label{units1}
The set ${\cal V}_0$ of elements of norm 1 of ${\cal Q}_0$ of the form 
\[ [r(x,z)+s(x,z){\bf i}+t(x,z){\bf j}+w(x,z){\bf k}],\]where
\begin{equation}
\label{integercoefficients }
2r,2s,2t,2w \in \IZ[x,z],
\end{equation}
\begin{equation}
\label{r00}
r(0,0) =1
\end{equation}
and
\begin{equation}
\label{oldcongruence}
s(x,z) \equiv  zw(x,z) \;\; \mbox{mod} \; x,
\end{equation}
is a group.  This group is not trivial as,  for instance,
\[ \Big\{ \frac{1}{2}\Big(x+2+x{\bf i}\Big),   \frac{1}{2}\Big(x+2+(x-2z){\bf i}-2{\bf k}\Big),  \frac{1}{2}\Big(z+2 -z{\bf i}-{\bf j}-{\bf k}\Big) \Big\}\subset {\cal V}_0 \]
\end{theorem}

\noindent {\bf Proof.} Let ${\bf u}=(r,s,t,w)=(1/2)(r_1,s_1,t_1,w_1) \in {\cal V}_0$. The fact that this has norm 1 gives
${{\bf u}}^{-1}=\overline{{\bf u}}=(r,-s,-t,-w)$, so clearly ${\cal V}_0$ is closed under inversion. We need only show that it is also closed under multiplication.
We have
\begin{equation}
\label{det}
r_1^2-\left(\frac{x+4}{x}\right) s_1^2-z(z-x) t_1^2 +\left(\frac{x+4}{x}\right)z(z-x)w_1^2=4
\end{equation}
Reducing modulo 2
gives
\begin{eqnarray*}
{(r_1-s_1)}^2 &\equiv& r_1^2 -s_1^2 \\
&\equiv& z(z-x) (w_1^2-t_1^2) \\
&\equiv& z(z-x) {(w_1-t_1)}^2\;\;\;\;\;{\rm mod}\; 2,
\end{eqnarray*}
whence $r_1-s_1 \equiv w_1-t_1 \equiv 0 \;\;{\rm mod}\; 2$, and so
\begin{equation}
\label{oldcommon2}
r_1(x,z) \equiv s_1(x,z) \;\; \mbox{mod} \; 2,
\end{equation}
\begin{equation}
\label{oldcommon22}
t_1(x,z) \equiv w_1(x,z) \;\; \mbox{mod} \; 2.
\end{equation}
Let ${\bf u_2}=(1/2)(r_2,s_2,t_2,w_2),{\bf u_3}=(1/2)(r_3,s_3,t_3,w_3) \in U$, then using (\ref{quatmult}), ${\bf u_2}{\bf u_3}=(1/2)(r,s,t,w)$, where
\begin{eqnarray}
2r&=&r_2r_3+\left(\frac{x+4}{x}\right)s_2s_3+z(z-x)t_2t_3-z(z-x)\left(\frac{x+4}{x}\right)w_2w_3 \nonumber \\
&=&r_2r_3+z(z-x)t_2t_3+(x+4)z w_2 w_3 \nonumber \\
\label{2r}
& & + (x+4)\left( \frac{s_2(s_3-zw_3)}{x}+ \frac{zw_3(s_2-zw_2)}{x} \right)  \\
2s&=&r_2s_3+s_2r_3+z(z-x)(w_2t_3-t_2w_3) \nonumber \\
2t&=&r_2t_3+t_2r_3+\left(\frac{x+4}{x}\right)(s_2w_3-w_2s_3) \nonumber \\
2w&=&r_2w_3+w_2r_3+s_2t_3-t_2s_3 \nonumber
\end{eqnarray}
The congruence (\ref{oldcongruence}) applied to ${\bf u_2}$ and ${\bf u_3}$ shows that each of these is a polynomial.  Moreover, since (\ref{oldcongruence}) also gives $s_2(0,0)=0$, setting $x=z=0$ in (\ref{2r}) then gives (\ref{r00}). Next 
\[
2(s-zw)=(r_2+zt_2)(s_3-zw_3)+(r_3-zt_3)(s_2-zw_2)+xz(t_2w_3-w_2t_3),
\]
so (\ref{oldcongruence}) holds for ${\bf u_2}{\bf u_3}$. It remains to show that the polynomials $r$, $s$, $t$ and $w$ have integer coefficients. We have
\begin{eqnarray}
2r&=&r_2r_3+\left(\frac{x+4}{x}\right)s_2s_3+z(z-x)t_2t_3-z(z-x)\left(\frac{x+4}{x}\right)w_2w_3 \nonumber \\
&\equiv& r_2r_3+s_2s_3+z(z-x)(t_2t_3-w_2w_3) \;\; \mbox{mod} \; 2 \nonumber \\
&\equiv& 0\;\; \mbox{mod} \; 2,
\end{eqnarray}
using (\ref{oldcommon2}) and (\ref{oldcommon22}), whence $r \in \IZ[x,z]$. Similar, and easier, arguments give the same conclusion for $s$, $t$ and $w$. \hfill $\Box$

\bigskip

\noindent {\bf Remark.} It is not difficult to show that
$r(0,0) = \pm 1$ and $s(x,z) \equiv  \pm zw(x,z) \;\; \mbox{mod} \; x$ follow from (\ref{det}), so that (\ref{r00}) and
(\ref{oldcongruence}) are just normalizing choices of sign.

\bigskip

For each fixed $\beta, \beta', \gamma \in \
\IC$, with $\beta \neq 0$, if we let $Q=\sqrt{\beta(\beta+4)}$ and $D_1$ and $D_2$ be any fixed numbers such that $D_1D_2=\gamma(\gamma-\beta)/\beta^2$,
then the evaluation map
\begin{eqnarray}
\lefteqn{ \phi_{\beta, \beta',\gamma}([r(x,z)+s(x,z){\bf i}+t(x,z){\bf j}+w(x,z){\bf k}])} \nonumber \;\;\;\;\;\;\;\;\;\; \\
\label{ev}
&=&\left[\begin{array}{cc}
r(\beta, \gamma)+\frac{s(\beta, \gamma)Q}{\beta} & D_1(\beta t(\beta, \gamma)+w(\beta, \gamma)Q)  \\ D_2(\beta t(\beta, \gamma)-w(\beta, \gamma)Q) & r(\beta, \gamma)-\frac{s(\beta, \gamma)Q}{\beta}
\end{array}
\right]
\end{eqnarray}
is an algebra homomorphism from ${\cal Q}_0$ to $M_2(\IC)$, the algebra of $2 \times 2$ matrices over $\IC$.

For $\beta=0$, $\gamma \neq 0$, we set
\begin{eqnarray}
\lefteqn{ \phi_{0,\beta' \gamma}([r(x,z)+s(x,z){\bf i}+t(x,z){\bf j}+w(x,z){\bf k}]) } \;\;\;\;\;\;\;\; \nonumber \\
\label{evzero}
&=&\left[\begin{array}{cc}
r(0, \gamma)+\gamma t(0, \gamma) & 4g(0, \gamma)+2w(0,\gamma)  \\ 2\gamma w(0,\gamma) & r(0, \gamma)-\gamma t(0, \gamma)
\end{array}
\right],
\end{eqnarray}
where $g$  (a polynomial by (\ref{oldcongruence})) is given by
\begin{equation}
\label{gdef}
 g(x,z):=\frac{s(x,z)-zw(x,z)}{x}, 
 \end{equation}
and
\begin{eqnarray}
\lefteqn{ \phi_{0,\beta', 0}([r(x,z)+s(x,z){\bf i}+t(x,z){\bf j}+w(x,z){\bf k}])}\;\;\;\;\;\;\;\; \nonumber \\
\label{evzerozero}
&=&\left[\begin{array}{cc}
r(0,0) & 4g(0,0)-\left(\beta'+\sqrt{\beta'}\sqrt{\beta'+4}\right)w(0,0)
 \\ 0 & r(0,0)
\end{array}
\right]
\end{eqnarray}

The next theorem shows that the maps $\phi_{0, \beta', \gamma}$ arise as limits of maps $\phi_{\beta, \beta', \gamma}$ (after conjugating $\phi_{\beta, \beta', \gamma}$ in such a way as to make its fixed points approach a common limit as $\beta \to 0$). It follows that the maps $\phi_{\beta, \beta', \gamma}$ are all algebra homomorphisms (this is also not difficult to show directly).
In particular each $\phi_{\beta, \beta', \gamma}$ restricted to ${\cal V}_0$ is a group homomorphism to $SL(2, \IC)$, and thence by projection to $PSL(2, \IC)$.

\begin{theorem}
\label{limit}
Suppose that $\beta\neq 0$,  $\gamma \neq 0$, $k^2=1/\sqrt{\beta}  $,  and that
\begin{eqnarray*}
 Q=\sqrt{\beta(\beta+4)}, && Q'=\sqrt{\beta}\sqrt{\beta+4},\\
m= -k^{-1}\left[\frac{1}{\sqrt{\beta}}+\frac{\sqrt{\beta'+4}}{2\sqrt{\gamma}}\right], &&
m_1= \frac{\sqrt{\beta'+4}+\sqrt{\beta'}}{2\sqrt{\gamma}},\\
D_1=-\frac{1}{2}\sqrt{\frac{\gamma}{\beta}}
\left(\sqrt{\beta'+4}+\sqrt{\frac{4 \gamma+\beta \beta'}{\beta}}\right),
&&
D_2=\frac{1}{2}\sqrt{\frac{\gamma}{\beta}}
\left(\sqrt{\beta'+4}-\sqrt{\frac{4 \gamma+\beta \beta'}{\beta}}\right)
\end{eqnarray*}
(i.e. $D_1=ab$ and $D_2=cd$, where $a$, $b$, $c$, $d$ are given by (\ref{ad}) and (\ref{gammanonzerooffdiagonals})),
\begin{eqnarray*}
D'_1& = & -\left(\frac{\sqrt{\gamma}}{2\sqrt{\beta}}\right)
\left(\sqrt{\beta'+4}+\frac{2\sqrt{\gamma}\sqrt{1+\beta \beta'/(4 \gamma)}}{\sqrt{\beta}}\right),\\
D'_2&=&\left(\frac{\sqrt{\gamma}}{2\sqrt{\beta}}\right)
\left(\sqrt{\beta'+4}-\frac{2\sqrt{\gamma}\sqrt{1+\beta \beta'/(4 \gamma)}}{\sqrt{\beta}}\right),
\\
C&=&
\left\{\begin{array}{ll} \left[\begin{array}{cc}
\sqrt{D'_1/D_1} & 0 \\ 0 & \sqrt{D'_2/D_2}
\end{array}
\right], \;\;\;\;\;\; \mbox{if $Q'=Q$} \\ \\
\left[\begin{array}{cc}
0 & i \sqrt{D'_1/D_2}  \\ i \sqrt{D'_2/D_1} & 0
\end{array}
\right], \;\;\;\mbox{if $Q' =- Q$}
\end{array}
\right.
\\
M=\left[\begin{array}{cc}
k & m \\ 0 & k^{-1}
\end{array}
\right]C, && M_1=\left[\begin{array}{cc}
1 & m_1 \\ 0 & 1
\end{array}
\right].
\end{eqnarray*}
Then for $\gamma \neq 0$, ${\bf x} \in {\cal Q}$,
\begin{equation}
\label{firstlimit}
\lim_{\beta \to 0} M\phi_{\beta, \beta', \gamma}({\bf x})M^{-1}=\phi_{0,\beta', \gamma}({\bf x})
\end{equation}
and
\begin{equation}
\label{secondlimit}
\lim_{\gamma \to 0} M_1\phi_{0, \beta', \gamma}({\bf x})M_1^{-1}=\phi_{0,\beta', 0}({\bf x})
\end{equation}
\end{theorem}

\noindent {\bf Proof.} Since $D'_1D'_2=D_1D_2$, the diagonal entries of $C$ are the same, in the case $Q'=Q$, and the off-diagonal entries of $C$ divided by $i$ are mutually reciprocal, otherwise. Thus we can apply (\ref{generalconjuation}) and (\ref{offdiag}) respectively to obtain
\begin{eqnarray*}
&&C \phi_{\beta, \beta' \gamma}({\bf x})
C^{-1} =\left[\begin{array}{cc}
r(\beta, \gamma)+\frac{s(\beta, \gamma)Q'}{\beta} & D'_1(\beta t(\beta, \gamma)+w(\beta, \gamma)Q')  \\ D'_2(\beta t(\beta, \gamma)-w(\beta, \gamma)Q') & r(\beta, \gamma)-\frac{s(\beta, \gamma)Q'}{\beta}
\end{array}
\right]
\end{eqnarray*}
Thus
\begin{equation}
M\phi_{\beta, \beta' \gamma}({\bf x})M^{-1}=\left[\begin{array}{cc}
a_{11} & a_{12} \\ a_{21} & a_{22}
\end{array}
\right],
\end{equation}
where, using (\ref{generalconjuation}), and writing $r(\beta, \gamma)=r$ etc.
\begin{eqnarray*}
\lefteqn{a_{11} = r+\frac{sQ'}{\beta}  -k^{-2}\left[\frac{1}{\sqrt{\beta}}+\frac{\sqrt{\beta'+4}}
{2\sqrt{\gamma}}\right]\left(\frac{\sqrt{\gamma}}{2\sqrt{\beta}}\right) \cdot }\\ &&  \cdot
\left(\sqrt{\beta'+4}-\frac{2{\sqrt{\gamma}}\sqrt{1+\beta \beta'/(4 \gamma)}}{\sqrt{\beta}}\right)(\beta t-wQ') \\
&=& r+\frac{sQ'}{\beta}-\frac{\sqrt{\gamma}}{2}\left[1+\frac{\sqrt{\beta} \sqrt{\beta'+4}}
{2\sqrt{\gamma}}\right] \cdot \\
&& \cdot \left(\sqrt{\beta'+4}-2\frac{\sqrt{\gamma}}{\sqrt{\beta}}+O(\sqrt{\beta})\right)\left(\sqrt{\beta} t-2w+O(\beta)\right)
\end{eqnarray*}
Thus
\begin{eqnarray*}
\lim_{\beta \to 0} a_{11} &=& \lim_{\beta \to 0} \left[ r+\frac{sQ'}{\beta}-\frac{\sqrt{\gamma}}{2}\left(1+\frac{\sqrt{\beta} \sqrt{\beta'+4}}
{2\sqrt{\gamma}}\right)
\left(\sqrt{\beta'+4}-2\frac{\sqrt{\gamma}}{\sqrt{\beta}}\right)(\sqrt{\beta} t-2w)\right] \\
&=& \lim_{\beta \to 0} \left[ r+\frac{2s}{\sqrt{\beta}}+\sqrt{\gamma}\left(1+\frac{\sqrt{\beta} \sqrt{\beta'+4}}
{2\sqrt{\gamma}}\right)
\left(t\sqrt{\gamma}+w\sqrt{\beta'+4}-
2w\frac{\sqrt{\gamma}}{\sqrt{\beta}}
\right)\right]\\
&=& \lim_{\beta \to 0} \left[ r+\frac{2s}{\sqrt{\beta}}+\sqrt{\gamma}\left(
t\sqrt{\gamma}
+w\sqrt{\beta'+4}
-w\sqrt{\beta'+4}
\right)-\frac{2w \gamma}{\sqrt{\beta}} \right]\\
&=& \lim_{\beta \to 0} \left[ r+\frac{2(s-\gamma w)}{\sqrt{\beta}}+
t\gamma \right]=r(0,\gamma)+
\gamma t(0,\gamma),
\end{eqnarray*}
using (\ref{oldcongruence}) at the last step.
Since conjugation preserves traces, we then have $ \lim_{\beta \to 0} a_{22}=r(0,\gamma)-\gamma t(0,\gamma)$.
\begin{eqnarray*}
a_{21} &=& k^{-2}\left(\frac{\sqrt{\gamma}}{2\sqrt{\beta}}\right)
\left(\sqrt{\beta'+4}-\frac{2{\sqrt{\gamma}}\sqrt{1+\beta \beta'/(4 \gamma)}}{\sqrt{\beta}}\right)(\beta t-wQ') \\
&=& \left(\frac{\sqrt{\gamma}}{2}\right)
\left(\sqrt{\beta'+4}-\frac{2{\sqrt{\gamma}}\sqrt{1+\beta \beta'/(4 \gamma)}}{\sqrt{\beta}}\right)(\beta t-wQ')
\end{eqnarray*}
so
\begin{equation}
\lim_{\beta \to 0} a_{21} = \left(\frac{\sqrt{\gamma}}{2}\right) \lim_{\beta \to 0} \left(\frac{-2\sqrt{\gamma}}{\sqrt{\beta}}\right)
\left(-2w\sqrt{\beta}\right)=2\gamma w(0,\gamma)=2s(0,\gamma),
\end{equation}
again using (\ref{oldcongruence}).

Finally we show that $
\lim_{\beta \to 0} a_{12}=2w(0,\gamma)+4g(0,\gamma)
$.
Since ${\rm Det}(\phi_{\beta, \beta', \gamma})=1$, and this determinant is preserved under  conjugation and limits, it suffices to show that
\[
(r(0,\gamma)+\gamma t(0,\gamma))(r(0,\gamma)-\gamma t(0,\gamma))-2s(0,\gamma)(2w(0,\gamma)+4g(0,\gamma))=1.
\]
This is readily verified by letting $\beta \to 0$ in (\ref{det}), keeping in mind the definition of $g$, (\ref{gdef}).
This completes the proof of (\ref{firstlimit}). A similar, but much easier, calculation gives (\ref{secondlimit}).
Forming the conjugate
\[
\left[\begin{array}{cc}
1 & m_1  \\ 0 & 1
\end{array}
\right]
\left[\begin{array}{cc}
r(0, \gamma)+\gamma t(0, \gamma) & 4g(0,\gamma)+2w(0,\gamma)  \\ 2\gamma w(0,\gamma) & r(0, \gamma)-\gamma t(0, \gamma)
\end{array}
\right]
\left[\begin{array}{cc}
1 & -m_1  \\ 0 & 1
\end{array}
\right],
\]
using (\ref{generalconjuation}), and letting $\gamma \to 0$, gives the matrix at (\ref{evzerozero}). \hfill $\Box$

\bigskip

\subsection{ Proof of Theorem \ref{matrixrep}}
Let
\begin{eqnarray}
{\bf w_1}&=& \frac{1}{2}(x+2,x,0,0), \nonumber \\
{\bf w_2}&=& \frac{1}{2}(x+2,x-2z,0,-2) \nonumber\\
\label{ws}
{\bf w_3}&=& \frac{1}{2}(z+2,-z,-1,-1).
\end{eqnarray}
Each ${\bf w_i} \in {\cal V}_0$. Let $A$, $B$ be as in Theorem \ref{canonicalgroup}, $Q=\sqrt{\beta(\beta+4)}$. In the first case, $\beta \neq 0$, we calculate
\begin{equation}
\label{evenpower}
A^2=\frac{1}{2}\left[\begin{array}{cc}
\beta+2+Q & 0 \\0 & \beta+2-Q
\end{array}
\right]
\end{equation}

\begin{equation}
\label{conjevenpower2}
BA^2B^{-1}=\frac{1}{2}\left[\begin{array}{cc}
\beta+2+(\beta-2\gamma)Q /\beta & -2abQ \\ 2cdQ & \beta+2-(\beta-2\gamma)Q /\beta
\end{array}
\right]
\end{equation}

\begin{eqnarray}
\label{commutator}
[B,A]&=&\frac
{1}{2}\left[\begin{array}{cc}
\gamma+2-\gamma Q / \beta
& -ab(\beta+Q)
\\ -cd(\beta-Q) & \gamma+2+\gamma Q / \beta,
\end{array}
\right].
\end{eqnarray}
In the second case, $\beta=0$, $\gamma \neq 0$, we have
\begin{equation}
\label{evenpowerpar}
A^2=
\left[\begin{array}{cc}
1 & 2 \\0 & 1
\end{array}
\right],
\end{equation}

\begin{equation}
\label{conjevenpowerpar}
BA^2B^{-1}=\left[\begin{array}{cc}
1 & 0 \\ -2\gamma & 1
\end{array}
\right],
\end{equation}

\begin{equation}
\label{commutatorpar}
[B,A]=
\left[ \begin{array}{cc} 1 & -1 \\ -\gamma & \gamma+1 \end{array} \right]
\end{equation}

Finally, if $\beta=\gamma=0$, then $A^2$ is still given by (\ref{evenpowerpar}), and
\begin{equation}
\label{conjevenpowerparspecial}
BA^2B^{-1}=\left[\begin{array}{cc}
1 & 2+\beta(g)+\sqrt{\beta(g)}\sqrt{\beta(g)+4} \\ 0 & 1
\end{array}
\right] 
\end{equation}
and
\begin{equation}
[B,A]=\left[\begin{array}{cc}
1 & \left(\beta(g)+\sqrt{\beta(g)}\sqrt{\beta(g)+4}\right)/2 \\ 0 & 1
\end{array}
\right]
\end{equation}
In this last case the matrices do depend on $\beta(g)$, but are independent of the parameter $\ell$.

Thus, in all cases, a straightforward calculation using the evaluation map at (\ref{ev}) gives
\[
A^2=\phi_{\beta,\beta',\gamma}({\bf w}_1),\;\;\;BA^2B^{-1}=\phi_{\beta,\beta',\gamma}({\bf w}_2),\;\;\;[B,A]=\phi_{\beta,\beta',\gamma}({\bf w}_3),
\]
where $\beta=\beta(f)$, $\beta'=\beta(g)$ and $\gamma=\gamma(f,g)$, and
where (in the case $\beta \neq 0$), we set $D_1=ab$ and $D_2=cd$ with $a,b,c,d$ given by (\ref{ad}) and (\ref{offdiagonals}).
Since, by Lemma \ref{evengen}, these words generate all regular balanced even words in $A$ and $B$, it follows immediately that every such word is $\phi_{\beta,\beta',\gamma}({\bf w})$ for some ${\bf w} \in {\cal V}_0$. Theorem \ref{matrixrep} then follows, using Theorem \ref{units1}, (\ref{oldcommon2}) and (\ref{oldcommon22}). \hfill $\Box$

\begin{corollary}
\label{limitingcases}
Let $A(\beta)$ and $B(\beta,\beta',\gamma)$ be the matrix representatives of $f$ and $g$ respectively given by Theorem \ref{canonicalgroup}, where now we have made the dependency on parameters $\beta=\beta(f)$, $\beta'=\beta(g)$ and $\gamma=\gamma(f,g)$ explicit.
Let $w(a,b)$ be a regular balanced even word. Let $\beta \neq 0$, and $k$, $M$ and $M_1$ as in Theorem \ref{limit},
then for $\gamma \neq 0$,
\begin{equation}
\label{firstlimit2}
\lim_{\beta \to 0} M(w(A(\beta), B(\beta,\beta',\gamma)))M^{-1}=w(A(0), B(0,\beta',\gamma))
\end{equation}
and
\begin{equation}
\label{secondlimit2}
\lim_{\gamma \to 0} M_1(w(A(0), B(0,\beta',\gamma)))M_1^{-1}=w(A(0), B(0,\beta',0))
\end{equation}
\end{corollary}
Since polynomials and the trace function are continuous, and trace is preserved under conjugation, Theorem
\ref{tracepolynomials} for $\beta=0$ follows from the case $\beta \neq 0$, by letting $ \beta \to 0$.

\subsection{A Change of Variable }
We now introduce two new parameters which can be used to describe  2-generator groups (up to conjugacy), and which, when $\beta \neq 0$, can be used interchangeably with $\beta$ and $\gamma$ and will simplify many formulas in what follows.
For $f,g \in Isom^+(\IH)$ we  define
\begin{equation}
\label{lambda}
\lambda=\lambda (f)=({\rm tr}^2(f)-2)/2=\cosh(\tau +i\eta),
\end{equation}
so $\lambda= (\beta(f)+2)/2$.  When $f$ is elliptic or loxodromic,
\begin{equation}
\label{mu}
\mu=\mu(f,g)=\frac{{\rm tr}^2(f)-2{\rm tr}[f,g]}{{\rm tr}^2(f)-4}.
\end{equation}
In terms of our earlier parameters
\[ \mu =1-\frac{2\gamma(f,g)}{\beta(f)} \]
Rewriting (\ref{betageom}) and (\ref{gammageom}) in terms of $\lambda(f)$ and $\mu(f,g)$ gives
\[
\lambda (f)=\cosh(\tau(f) +i\eta(f)),
\]
and
\[
\mu(f,g)=1-(\lambda(g)-1)\sinh^2(\Delta),
\]
where we recall $\Delta$ -- the complex distance between axes -- is defined above at (\ref{gammageom}).
An important special case is captured by the next lemma.
\begin{lemma}
If $g$ is order 2, $\lambda(g)=-1$, and we obtain the particularly simple form:
\[ \mu(f,g)=\cosh(2\Delta).\]
\end{lemma}

Unwinding these parameters gives
\begin{equation} 
\label{oldtonew}
\beta(f)=2(\lambda(f)-1), \;\;\;{\rm and}\;\;\; \gamma(f,g)=-(\lambda(f)-1)(\mu(f,g)-1) \end{equation}

 If $\lambda=\lambda(f)\neq 1$ then it determines $f$ up to congugacy, and if further, $\mu=\mu(f,g) \neq 1$, then $\lambda$ and $\mu$ together determine the group $\langle f,g \rangle$. When $\lambda \neq 1$, we can rewrite the matrices $A$ and $B$ at (\ref{AB}) in terms of the new parameters as
\begin{equation}
\label{newAB}
A = \left[\begin{array}{cc}
\frac{\sqrt{\lambda^2-1}+(\lambda-1)}{\sqrt{2(\lambda-1)}}
 & 0 \\ 0 & \frac{\sqrt{\lambda^2-1}-(\lambda-1)}{\sqrt{2(\lambda-1)}}
\end{array}
\right],\;\;  B=\left[\begin{array}{cc}
a & b \\ c & d
\end{array}
\right],
\end{equation}
where now, writing $\lambda(g)=\lambda'$,
\begin{equation}
\label{newad}
a=\frac{1}{\sqrt{2}}\left(\sqrt{\lambda'+1}+\sqrt{\lambda'-\mu}\right),
\;\;\;\;
d=\frac{1}{\sqrt{2}}\left(\sqrt{\lambda'+1}-\sqrt{\lambda'-\mu}\right),
\end{equation}
and $b$ and $c$ are given by (\ref{offdiagonals}) when $\mu=1$,
and otherwise
\begin{equation}
\label{newoffdiagonals}
c=-b= \sqrt{\frac{1-\mu}{2}}
\end{equation}
for $\mu \neq 1$.

\bigskip

Given a regular balanced even word $w$, we can now rewrite the matrix $w(A,B)$ at (\ref{matrixform1}) as
\begin{equation}
\label{matrixform2} \hspace{-1.5 em}
\left[\begin{array}{cc}
R_w(\lambda, \mu)+S_w(\lambda, \mu)\sqrt{\lambda^2-1} & 2ab(T_w(\lambda, \mu)+W_w(\lambda, \mu)\sqrt{\lambda^2-1})  \\ 2cd(T_w(\lambda, \mu)-W_w(\lambda, \mu)\sqrt{\lambda^2-1} & R_w(\lambda, \mu)-S_w(\lambda, \mu)\sqrt{\lambda^2-1}
\end{array}
\right]
\end{equation}
where $a$, $b$, $c$, $d$ are given by (\ref{newad}) and (\ref{newoffdiagonals}), and, setting $x=2(u-1)$, $z=-(u-1)(v-1)$, the polynomials $R_w, S_w, T_w, W_w$ are given by
\begin{equation}
\label{polyconvert1}
R_w(u,v)=r_w(x,z),  \;\;\; W_w(u,v)=w_w(x,z)
\end{equation}
\begin{equation}
\label{polyconvert2}
S_w(u,v)=2s_w(x,z)/x,  \;\;\; T_w(u,v)=xt_w(x,z)/2.
\end{equation}
The congruence (\ref{oldcongruence}) insures that $S_w$ is a polynomial. For arbitrary balanced words we take the above as {\em definitions} of $R_w$ etc. Recall that in this case there is a regular even balanced word $v$ such that $r_w=r_v$, $s_w=s_v$ etc., so that it remains true that $S_w$ is a polynomial in this case.

\medskip

Each of the polynomials $2R_w$, $2S_w$, $2T_w$, and $2W_w$ has integer coefficients,  and
corresponding to these new parameters, we define the quaternion algebra ${\cal Q}$ by
\begin{equation}\label{quatQ}
{\cal Q}=\left( \frac{ u^2-1, v^2-1}{\IR(u,v)} \right),
\end{equation}
Here, the indeterminates  $u$ and $v$ correspond to $\lambda$ and $\mu$ respectively.
It is straightforward to show that the map $\rho : {\cal Q}_0 \to {\cal Q}$ given by
\begin{equation}
\label{outconvert}
\rho(r ,s ,t ,w ) = (r ,\frac{s }{u-1},(u-1)t , w )
\end{equation}
is an isomorphism.
On the right hand side $x$ and $z$ are converted into terms of $u$ and $v$ by the formulae
\begin{equation}
\label{conversions}
x=2(u-1), \;\; z=-(u-1)(v-1), \;\;\;  \left( u=\frac{x+2}{2}, \;\;\; v=1-\frac{2z}{x} \right),
\end{equation}
these conversions being just the same as those relating $\beta$ and $\gamma$ to $\lambda$ and $\mu$.
The proof uses the observations that
\[
u^2-1=\frac{x^2}{2} \left(\frac{x+4}{x}\right)\;\;\mbox{and}\;\;  v^2-1=\frac{2}{x^2}\; z(z-x) .
\]
It is then easy to see that the inverse map is given by
\begin{equation}
\label{inverse}
\rho^{-1}(R, S , T , W ) = (R , \frac{x}{2}S , \frac{2}{x}T , W )
\end{equation}
Where we now use the second pair of equations in (\ref{conversions}) to convert the right hand side back into terms of
$x$ and $z$.

\bigskip
 We can now characterize the image under $\rho$ of the group ${\cal V}_0$ defined in Theorem \ref{units1}. This result is a direct consequence of (\ref{outconvert}), (\ref{conversions}), (\ref{inverse}) and the definition of ${\cal V}_0$.
\begin{theorem}
\label{units2}
${\cal V}:=\rho({\cal V}_0)$ comprises the elements
\[ (R,S,T,W) = (R(u,v), S(u,v), T(u,v), W(u,v)) \]
 of ${\cal Q}$ for which $R(1,1)=1$, and each of $2R$, $2(u-1)S$, $2{(u-1)}^{-1}T$, $2W$ and $S+(v-1)W$ is a polynomial of the form
\[
\sum a_{n,m} {(u-1)}^m {(v-1)}^n
\]
such that, for each term, $m \geq n$ and $a_{n,m}$ is an integer multiple of $2^{m-n}$ (in particular, each $a_{n,m}$ is an integer).
\end{theorem}

\noindent{\bf Remark.} The condition on $S(u,v)+(v-1)W(u,v)$ is equivalent to (\ref{oldcongruence}): we have $g(x,z):=(s(x,z)-zw(x,z))/x=(S(u,v)+(v-1)W(u,v))/2$. Note also that this condition insures that $S(u,v)$ is a polynomial.

\bigskip

If ${\bf x}=(r(x,z),s(x,z),t(x,z),w(x,z)) \in {\cal V}_0$, and if
\[ \rho({\bf x})=(R(u,v),S(u,v),T(u,v),W(u,v)) \in {\cal V},\]
then by the definition of $\rho$, the matrix at (\ref{matrixform2}) is $\phi_{\beta, \beta',\gamma}({\bf x})$ ($\beta \neq 0$).
Accordingly we define,
for each $\lambda \neq 1,\lambda',\mu \in \IC$, the algebra homomorphisms $\psi_{\lambda, \lambda',\mu}:  {\cal Q} \to M_2(\IC)$ by
\begin{eqnarray}
\lefteqn{ \psi_{\lambda, \lambda',\mu}(R,S,T,W)}&& \label{newev} \\ \nonumber\\ \nonumber
\label{newev}
&=&\left[\begin{array}{cc}
R(\lambda, \mu)+S(\lambda, \mu)\sqrt{\lambda^2-1} & 2ab(T(\lambda, \mu)+W(\lambda, \mu)\sqrt{\lambda^2-1})  \\
2cd(T(\lambda, \mu)+W(\lambda, \mu)\sqrt{\lambda^2-1}) & R(\lambda, \mu)-S(\lambda, \mu)\sqrt{\lambda^2-1}
\end{array}
\right], \nonumber
\end{eqnarray}
where $a$, $b$, $c$, $d$ are given by (\ref{newad}) and (\ref{newoffdiagonals}). We thus have
\begin{proposition}
For $\beta \neq 0, \beta', \gamma \in \IC$, $x \in {\cal V}_0$,
\[
\psi_{\lambda, \lambda',\mu}(\rho({\bf x}))=\phi_{\beta, \beta',\gamma}({\bf x}),
\]
where $\lambda=1+\beta/2$,
$\lambda'=1+\beta'/2$ and $\mu=1-2\gamma/\beta$.
\end{proposition}
\bigskip

In some respects $\lambda$ and $\mu$ are better parameters to use than $\beta$ and $\gamma$: they have a simpler geometrical interpretation, the matrix representations and quaternion algebras are simpler and neater, and there is an
obvious symmetry between $\lambda$ and $\mu$, corresponding to the symmetry between two loxodromics with perpendicular axes (subsection \ref{asisometries} below). They also have a major drawback: $\mu$ is undefined when $f$ is parabolic, so to deal with this case we still need $\beta$ and $\gamma$.

\section{Elements of unit norm in ${\cal Q}$}
We have found a group ${\cal V}$ of elements of norm 1 in ${\cal Q}$, which maps under each evaluation homomorphism
$\psi_{\lambda, \lambda',\mu}$ to a group which includes the regular balanced even words in two
generators $f$ and $g$, where $f$ is elliptic or loxodromic.
In this section will show how ${\cal V}$ naturally extends to a larger group which we will denote ${\cal U}$, and this gives a corresponding extension of the isometry group $\psi_{\lambda, \lambda',\mu}({\cal V})$.
We look further at this group in \ref{asisometries}. To begin with we consider properties of elements of norm 1 in general.

The requirement that $(R,S,T,W) \in {\cal Q}$ has norm 1 is given expliciltly by
\begin{equation}
\label{symmdet}
R^2-(u^2-1)S^2-(v^2-1) T^2 +(u^2-1)(v^2-1)W^2=1
\end{equation}
We will confine our attention to the solutions of (\ref{symmdet}) for which $R,S,T,W$ are all polynomials, with the additional normalizing condition that $R(1,1)=1$. These solutions clearly form a group, which we denote by ${\cal U}_1$. 

\begin{table}[h]
\caption{$(R,S,T,W) \in {\cal U}_1$ of degree at most 2}
\centering
\begin{tabular}{|c|c|c|c|c|}
\hline
 $R$ & $S$ & $T$ & $W$ \\ \hline \hline
 1 & 0 & 0 & 0 \\ \hline
 $u$ & 1 & 0 & 0 \\ \hline
 $v$ & 0 & 1 & 0 \\ \hline
 $u$ & $v$ & 0 & 1\\ \hline
 $v$ & 0 & $u$ & 1\\ \hline
 $u v$ & 1 & $u$ & 0 \\ \hline
 $u v$ & $v$ & 1 & 0 \\ \hline
 $u v$ & $v$ & $u$ & 1\\ \hline
 $2 u^2-1$ & $2 u$ & 0 & 0 \\ \hline
 $ 2 v^2-1$ & 0 & $2v$ & 0 \\ \hline
 $(1 + u + v - u v)/2$ & $(v-1)/2$ & $(u-1)/2$ & $1/2$ \\ \hline
 $(1 + u - v + u v)/2$ & $(v+1)/2$ & $(u-1)/2$ & $1/2$\\ \hline
 $(1-u+v+uv)/2$ & $(v-1)/2$ & $(u+1)/2$ & $1/2$\\ \hline
 $(-1 + u + v + u v)/2$ & $(v+1)/2$ & $(u+1)/2$ & $1/2$\\ \hline
\end{tabular}
\label{table:quadraticunits}
\end{table}


\medskip

We define the {\em degree} of ${\bf u}=(R,S,T,W)\in {\cal U}_1$ by
\[ {\rm deg}({\bf u})=\max \{{\rm deg}(R),{\rm deg}(S)+1,{\rm deg}(T)+1, {\rm deg}(W)+2 \}.\]
 It is easy to check that ${\rm deg}({\bf u}{\bf v}) \leq {\rm deg}({\bf u})+{\rm deg}({\bf v})$. 
 For any fixed degree it is possible in principle to evaluate all members of ${\cal U}_1$ of any fixed degree $d$, by equating coefficients in (\ref{symmdet}), and we have done this for $d \leq 4$. 
 Table \ref{table:quadraticunits} lists all members of ${\cal U}_1$ of degree at most 2 (up to sign changes of the components $S$, $T$ and $W$). In this case, the polynomials, like those of ${\cal V}$, all have integer or half-integer coefficients. Our main result in this section is that 
all members of ${\cal Q}$ with components in $\IQ[u,v]$ and  norm in $\IZ[u,v]$ have this property, and that they form an order (with reference to the underlying ring $\IZ[u,v]$).
We define
\[
{\cal O}=\{{\bf u}=(R,S,T,W) \in {\cal Q}\;|\; R,S,T,W \in \IQ[u,v], N({\bf u}) \in \IZ[u,v] \}
\]
\begin{theorem}
\label{halfintgroup}
\begin{enumerate}
\item ${\cal O}$ is the set of quaternions of the form $(R,S,T,W)+\frac{1}{2}P((u+1)(v+1),v+1,u+1,1)$,
where $R,S,T,W,P \in \IZ[u,v]$.
\item ${\cal O}$ is the unique maximal order of ${\cal Q}$, which contains ${\bf i}$ and ${\bf j}$.
\end{enumerate}
\end{theorem}

It follows (since ${\cal O}$ is clearly closed under conjugation) that the elements of
${\cal O}$ of norm 1 form a group, which we denote ${\cal U}$.

The appearance of half-integer coefficients in ${\cal O}$ is reminiscent of the {\em Hurwitz order} ${\cal H}$ in $\IHM$, the quaternions of Hamilton, defined as ${\cal H}=\{\frac{1}{2}(n_1+n_2{\bf i}+n_3{\bf j}+n_4{\bf k}) \in \IHM\;|\; 
n_1,n_2,n_3,n_4 \in \IZ, n_1 \equiv n_2 \equiv n_3 \equiv n_4 \bmod 2\}$. The definition of ${\cal O}$ has no analog for ${\cal H}$; there are plenty of quaternions with integer norm and  components which are rational but not half-integers, for example
$(3/5,4/5,0,0)$.
We do however have the following counterpart to Theorem \ref{halfintgroup} (2), see e.g. \cite{V}

\begin{theorem}
\label{Hurwitzorder}
${\cal H}$ is the unique maximal order of $\IHM$, which contains ${\bf i}$ and ${\bf j}$.
\end{theorem}
We will take an axiomatic approach which covers both of these theorems.

We may characterize ${\cal H}$ as the set of quaternions with integer norm and half-integer components. This amounts to the simple observation that, for integers $a,b,c,d$ 
\begin{equation}
\label{foursquares}
a^2+b^2+c^2+d^2 \equiv 0 \bmod 4 \Rightarrow a \equiv b \equiv c \equiv d \bmod 2.
\end{equation}
The integers also satisfy the similar property
\begin{equation}
\label{strongfoursquares}
a^2+b^2+c^2+d^2 \equiv 0 \bmod 8 \Rightarrow a \equiv b \equiv c \equiv d \equiv 0\bmod 2.
\end{equation}
As a simple application of this we observe that a quaternion ${\bf u}$ with rational components and integer norm has no component with denominator divisible by 4. For if this occurred we would have, clearing denominators, a quaternion with integer components not all even and norm divisible by 16, contrary to (\ref{strongfoursquares}).

\medskip

We say that a commutative ring which satisfies (\ref{foursquares}) or (\ref{strongfoursquares}) has respectively the
{\em four squares} property and the {\em strong four squares} property. To justify this terminology, we show that $(\ref{strongfoursquares}) \Rightarrow (\ref{foursquares})$.  Suppose that (\ref{strongfoursquares}) holds, and that 
 $a^2+b^2+c^2+d^2 \equiv 0 \bmod 4$, then ${(a-b)}^2+{(a+b)}^2+{(c-d)}^2+{(c+d)}^2=2(a^2+b^2+c^2+d^2) \equiv 0 \bmod 8$, and so applying (\ref{strongfoursquares}), we get $a \equiv b \bmod \;2$ and $c \equiv d \bmod \;2$. The same argument with $b$ and $c$ interchanged gives
 $a \equiv c \bmod \;2$, and so $a \equiv b \equiv c \equiv d \bmod \;2$, proving (\ref{foursquares}).
(The converse fails; consider for example $R=\IZ_4$.)

The main step in our proofs is to show that, if $R$ has the (strong) four squares property, then the polynomial rings $R[x_1,x_2, \ldots, x_n]$ also have this property, together with some generalisations  thereof.

\begin{lemma}
\label{sum4squaresetc}
Suppose $R$ is a commutative ring  and suppose $\varphi_i$ ($1 \leq i \leq 4$) are fixed polynomials in $R[x_1,x_2, \ldots, x_n]$ such that

\begin{enumerate}
\item the constant term in each $\varphi_i$ is $1$,
\item for each non-constant monomial $x_1^{p_1}x_2^{p_2} \ldots x_n^{p_n}$ in each $\varphi_i$, at least one of the powers $p_i$ is odd,
\item $\varphi_1 \equiv \varphi_2 \equiv \varphi_3 \equiv \varphi_4 \equiv 0 \bmod 2$,
\end{enumerate}
Then for all $k \geq 1$, if 
\begin{eqnarray}
\label{k1}
&& \mbox{For all $a \in R$, $a^2 \equiv 0 \bmod 2 \Rightarrow a \equiv 0 \bmod 2$, when $k=1$} \\
\label{keven}
&&\mbox{$R$  has the four squares property, when $k=2$} \\
\label{kodd}
&&\mbox{$R$  has the strong four squares property, when $k \geq 3$}
\end{eqnarray}
and, for $p_1,p_2,p_3,p_4 \in R[x_1,x_2, \ldots, x_n]$
\[ 
\varphi_1 p_1^2+\varphi_2 p_2^2+\varphi_3 p_3^2+\varphi_4 p_4^2 \equiv 0 \bmod 2^k,
\]
then
\begin{eqnarray*}
&\cdot &\mbox{$p_1+p_2+p_3+p_4 \equiv 0  \bmod 2$ when $k=1$} \\
&\cdot &\mbox{$p_1 \equiv p_2 \equiv p_3 \equiv p_4 \bmod 2^{k/2}$ when $k$ is even} \\
&\cdot  &\mbox{$p_1 \equiv p_2 \equiv p_3 \equiv p_4 \equiv 0 \bmod 2^{(k-1)/2}$ when $k$ is odd}
\end{eqnarray*}
\end{lemma}

\smallskip

\noindent {\bf Proof.} 
First we note that if $R$ has a four squares property, then (\ref{k1}) holds, for
 if (\ref{k1}) fails then there is $r \in R$ with $r^2 \equiv 0 \bmod 2$, $r \not\equiv 0 \bmod 2$, in which case (\ref{keven}) fails with $a=b=r$, $c=d=0$. 

 For convenience, we suppose that $n=2$ (the proof for $n>2$ is an obvious generalization of this). Throughout this proof we order $\IZ^2$ lexicographically, that is $(a,b) < (c,d)$ when either $a<c$ or $a=c$ and $b<d$.

Suppose that polynomials $\varphi_i(x,y)=\sum c_i(m,n)x^m y^n \in R[x,y]$ are as in Lemma \ref{sum4squaresetc}, and that $p_i(x,y)=\sum a_i(m,n)x^m y^n \in R[x,y]$.

We use induction on $k$. Let $k \in \IN$, and suppose that the theorem holds for smaller values.
The hypotheses are
\begin{equation}
\label{full4squares}
\varphi_1 (x,y) p_1^2(x)+\varphi_2 (x,y)p_2^2(x)+\varphi_3 (x,y)p_3^2(x)+\varphi_4 (x,y)p_4^2(x) \equiv 0 \bmod 2^k.
\end{equation}
together with the conditions (\ref{k1}), (\ref{keven}) and (\ref{kodd}) on $R$ according as $k=1$, $k=2$ or $k \geq 3$. First suppose that $k=k_1>3$, then by the case $k=3$ each $p_i \equiv 0 \bmod \;2$, and by the induction hypothesis we apply the case $k=k_1-2$ to the $p_i/2\in R[x,y]$ to get the required result.
We suppose then that $k \leq 3$.
 
 For $k=1,2,3$ respectively, the required result can be stated in terms of coefficients as, for all $n,m \in \IZ$,
\begin{eqnarray}
\label{equalcoeffs1}
&\cdot &a_1(n,m)+a_2(n,m)+a_3(n,m)+ a_4(n,m) \equiv 0 \bmod 2 \\ 
\label{equalcoeffs}
&\cdot &a_1(n,m) \equiv a_2(n,m)\equiv a_3(n,m)\equiv a_4(n,m) \bmod 2 \\ 
\label{equalcoeffs3}
&\cdot &a_1(n,m) \equiv a_2(n,m)\equiv a_3(n,m)\equiv a_4(n,m) \equiv 0 \bmod 2.
\end{eqnarray}

To make the induction go through we will prove in the case $k=3$ that, in addition to
(\ref{equalcoeffs3}), 
\begin{equation}
\label{sumcoeffs}
a_1(n,m)+a_2(n,m)+ a_3(n,m)+a_4(n,m) \equiv 0 \bmod 4 
\end{equation}

We set $p^2_i(x,y)=\sum s_i(m,n)x^m y^n$ ($1 \leq i \leq 4$), and define vectors
\begin{eqnarray*}
{\bf c}(m,n)&=&(c_1(m,n),c_2(m,n),c_3(m,n),c_4(m,n))\\
{\bf a}(m,n)&=&(a_1(m,n),a_2(m,n),a_3(m,n),a_4(m,n)) \\
{\bf s}(m,n)&=&(s_1(m,n),s_2(m,n),s_3(m,n),s_4(m,n)).
\end{eqnarray*}

We have, for $p$ and $q$ even,
\begin{eqnarray}
\label{evenpower}
s_i(p,q) &=& 2\sum_{\scriptsize{\begin{array}{c} (s,t) < (p/2, q/2)  \end{array} }} a_i(s,t)a_i(p-s,q-t)+ a_i^2(p/2,q/2)
\end{eqnarray}
and, for $p$ or $q$ odd
\begin{eqnarray}
\label{oddpower}
s_i(p,q) &=& 2\sum_{\scriptsize{\begin{array}{c} (s,t) < (p/2, q/2) \end{array} }}  a_i(s,t)a_i(p-s,q-t),
\end{eqnarray}
whereupon summing gives us
\begin{equation}
\label{squarecoefft}
\sum_{i=1}^4 s_i(p,q)=2\sum_{(s,t) < (p/2, q/2)} {\bf a}(s,t) \cdot {\bf a}(p-s,q-t)+\sum_{i=1}^4 a_i^2(p/2,q/2),
\end{equation}
with the last sum only present when $p$ and $q$ are even.

Equating the coefficient of $x^{2n} y^{2m}$ in the left side of (\ref{full4squares}) to 0 $\bmod \;     2^k$ gives, using the second hypotheses on the $\varphi_i$,
\begin{equation}
\label{coefft}
\sum_{\scriptsize{\begin{array}{c}
p+p'=2n,q+q'=2m \\ \mbox{$p$ or $q$ is odd}
\end{array}}}  {\bf s}(p,q) \cdot {\bf c}(p',q')+\sum_{i=1}^4 s_i(2n,2m) \equiv 0 \bmod 2^k.
\end{equation}
By (\ref{oddpower}), each term in the first sum is even, whence
using (\ref{squarecoefft}) with $p=2n$, $q=2m$, $\sum_{i=1}^4 a_i^2(n,m) \equiv 0 \bmod 2$. Since $\sum a_i^2(n,m) \equiv {\left( \sum a_i(n,m) \right)}^2 \bmod 2$, (\ref{k1}) gives the congruence (\ref{equalcoeffs1}). For $k=1$, this completes the proof.

\medskip

We now prove (\ref{equalcoeffs}) for $k=2$ and (\ref{equalcoeffs3}) and (\ref{sumcoeffs}) for $k=3$ by induction on $(n,m)$. Suppose the result holds for all $(s,t)<(n,m)$. We first show that 
\begin{eqnarray}
\label{sa0}
s_i(p,q) &\equiv& 0 \bmod \;2^{k-1},\;\; (i=1,2,3,4)\\
\label{sa1}
\sum_{i=1}^4 s_i(2n,2m) &\equiv& \sum_{i=1}^4 a_i^2(n,m) \bmod 2^k,  \;\;\; {\rm and} \\
\label{sa2}
\sum_{i=1}^4 s_i(p,q) &\equiv& 0 \bmod 2^k,
\end{eqnarray}
if $p \leq 2n$, $q \leq 2m$ and $p$ or $q$ is odd.  The first of these follows from (\ref{oddpower}) and the induction hypothesis, and the other two from (\ref{squarecoefft}) since, for $(s,t) < (p/2, q/2) \leq (n,m)$, and $k=2$
\begin{eqnarray*}
{\bf a}(s,t) \cdot {\bf a}(p-s,q-t) &\equiv &  a_1(s,t)\sum_{i=1}^4 a_i(p-s,q-t) \\
&\equiv& 0 \bmod 2,
\end{eqnarray*}
using the induction hypothesis at the first congruence, and (\ref{equalcoeffs1}) at the second.

For $k=3$, by the induction hypothesis, each $a_i(s,t)$ is even, for $(s,t) < (p/2, q/2) \leq (n,m)$, and
\begin{eqnarray*}
\frac{1}{2}{\bf a}(s,t) \cdot {\bf a}(p-s,q-t) &\equiv &  a_1(p-s,q-t)\sum_{i=1}^4 \frac{a_i(s,t)}{2} \\
&\equiv& 0 \bmod 2 \;\;\; \mbox{ by (\ref{sumcoeffs})}.
\end{eqnarray*}
Here we use the result for $k=2$, which gives (\ref{equalcoeffs}), at the first step and the induction hypothesis at the second.

Recalling that $c_1(p',q') \equiv 
c_2(p',q') \equiv c_3(p',q') \equiv c_4(p',q') \bmod 2$, by the third hypothesis on the $\varphi_i$, and using (\ref{sa0}), the summand in the first sum of (\ref{coefft}) is
\begin{eqnarray*}
{\bf s}(p,q) \cdot {\bf c}(p',q') &=& \sum_{i=1}^4 s_i(p,q) c_i(p',q') \\
&=&  2^{k-1}  \sum_{i=1}^4 \left( \frac{s_i(p,q)}{2^{k-1}}\right) c_i(p',q')\\
&\equiv&  2^{k-1}c_1(p',q') \sum_{i=1}^4 \frac{s_i(p,q)}{2^{k-1}} \bmod 2^k \\
&=&  c_1(p',q') \sum_{i=1}^4 s_i(p,q)\\
&\equiv& 0 \bmod 2^k \;\;\;\mbox{(by (\ref{sa2}))}
\end{eqnarray*}
This, together with (\ref{coefft}) and (\ref{sa1}) gives
$a^2_1(n,m)+ a^2_2(n,m)+ a^2_3(n,m)+ a^2_4(n,m) \equiv 0 \bmod 2^k$, whence by hypothesis
we get (\ref{equalcoeffs}) for $k=2$ and (\ref{equalcoeffs3})  for $k=3$,

When $k=3$ we have
${((a_1(n,m)+a_2(n,m)+a_3(n,m)+a_4(n,m))/2)}^2 \equiv
{(a_1(n,m)/2)}^2+{(a_1(n,m)/2)}^2+{(a_1(n,m)/2)}^2+{(a_1(n,m)/2)}^2 \equiv 0 \bmod 2$, whence (\ref{k1}) gives (\ref{sumcoeffs}). \hfill $\Box$
\medskip




\noindent  We first note a simple special case ($\varphi_1=\varphi_2=\varphi_3=\varphi_4=1$, $k=2,3$).
\begin{corollary}
\label{polynomialinheritance}
If the ring $R$ satisfies the (strong) four squares property (\ref {foursquares}), then so does the polynomial ring
$R [x_1, x_2, \ldots ,x_n]$. In particular this is true of $\IZ [x_1, x_2, \ldots ,x_n]$.
\end{corollary}

\smallskip


\begin{corollary}
\label{ax2congruences}
If $a,b \in \IZ[u_1,u_2 \ldots u_k]$ can be written 
$a \equiv \alpha  \alpha' \bmod 4$, $b \equiv \beta  \beta' \bmod 4$, with $\alpha,  \alpha', \beta,  \beta' \in \IZ[u_1,u_2 \ldots u_k]$ satisfying $\alpha \equiv \alpha' \bmod 2$, $\beta \equiv \beta' \bmod 2$, and $\varphi_1=\alpha \beta$, $\varphi_2=-\alpha' \beta$, $\varphi_3=-\alpha \beta'$, $\varphi_4=\alpha' \beta'$ satisfy the hypotheses 
of Lemma \ref{sum4squaresetc}, then for $R,S,T,W \in \IZ[u_1,u_2 \ldots u_k]$,
\begin{equation}
\label{4congr}
R^2-aS^2-bT^2+abW^2 \equiv 0 \bmod 4 \Rightarrow R \equiv aS \equiv bT \equiv ab W \bmod 2
\end{equation}
and if $a,b \not\equiv 0 \bmod 2$,
\begin{equation}
\label{8congr}
R^2-aS^2-bT^2+abW^2 \equiv 0 \bmod 8 \Rightarrow R \equiv S \equiv T \equiv W \equiv 0 \bmod 2
\end{equation}
\end{corollary}
\noindent {\bf Proof.} Multiplying the left hand side of (\ref{4congr}) through by $\alpha \beta$, and setting $r=R$, $s=\alpha S$, $t=\beta T$ and $w=\alpha \beta W$, gives
the equivalent form
\begin{equation}
\label{modified4norm}
\alpha \beta r^2-\alpha' \beta s^2-\alpha \beta' t^2+\alpha' \beta' w^2 \equiv 0 \bmod 4.
\end{equation}
Lemma \ref{sum4squaresetc}  with $k=2$ then gives (\ref{4congr}).
Similarly Lemma \ref{sum4squaresetc}  with $k=3$ gives (\ref{4congr}).
\hfill $\Box$

\bigskip

\begin{theorem}
\label{ax}
Let $R$ be an integral domain of characteristic $\neq 2$, with field of fractions $K$. Let $a,b, \alpha, \beta \in R$ be such that $\alpha^2 \equiv a \bmod 2$ and $\beta^2 \equiv b \bmod 2$, and let ${\cal Q}$ be the quaternion algebra
$\frac{a,b}{K}$. Let $O$ comprise the quaternions of the form ${\bf r}+\frac{r}{2}{\bf c}$, where ${\bf c}=(\alpha \beta, \beta, \alpha,1)$, and $r$ and the components of ${\bf r}$ are in $R$, then
\begin{enumerate}
\item $O$  is an order in ${\cal Q}$, and $N({\bf u}) \in R$ for each ${\bf u} \in O$. 
\item If further
\begin{enumerate}
\label{conditions}
\item  
\label{ic}
$R$ is integrally closed
\item 
\label{division}
2 is prime in $R$
\item 
\label{ab1}
$a$ and $b$ are not divisible by 2;
$a|x^2 \Rightarrow a|x$, $b|x^2 \Rightarrow b|x$ [whence
$a$ and $b$ are squarefree (i.e. have no square divisors other than units)].
\item 
\label{linind}
If $b|y^2-ax^2$, then $b|x,y$, and if
$a|y^2-bx^2$, then $a|x,y$.
\item 
\label{maincongruence}
If ${\bf u}=(x,y,z,w) \in R^4$, and $N({\bf u})=x^2-ay^2-bz^2+abw^2 \equiv 0 \bmod 4$, then
$x \equiv \alpha y \equiv \beta z \equiv \alpha \beta w \bmod 2$,
\end{enumerate}
then every quaternion in $\frac{1}{2}R^4$ with norm in $R$ is in $O$, and every order which contains ${\bf i}$ and ${\bf j}$ lies in $O$. In particular, $O$ is maximal.
\end{enumerate}
\end{theorem}

\bigskip

\noindent {\bf Proof.}  
Clearly $O$ is an ideal. For ${\bf u}=(x,y,z,w)  \in R^4$, a straightforward calculation gives ${\bf u}{\bf c} \equiv {\bf c}{\bf u} \equiv (x+\alpha y+\beta z+\alpha \beta w){\bf c} \bmod 2$ and 
${\bf c}^2=-(b-\beta^2)(a-\alpha^2)1+2\alpha \beta {\bf c} \equiv 2\alpha \beta {\bf c} \bmod 4$, from which it follows that $O$ is also a ring. Since $N({\bf c})=(\alpha^2-a)(\beta^2-b)\equiv 0 \bmod 4$, it readily follows that $N({\bf u}) \in R$ for ${\bf u} \in O$.

Now suppose that (\ref{ic})-(\ref{maincongruence}) hold. If ${\bf u}=(x,y,z,w)  \in \frac{1}{2}R^4$, then (\ref{maincongruence})
applied to $2{\bf u}$, gives $2x \equiv \alpha 2y \equiv \beta 2z \equiv \alpha \beta 2w \bmod 2$. By (\ref{division}) and (\ref{ab1}) we may cancel modulo 2 to obtain $2x \equiv 2\alpha \beta w$,
$2y \equiv 2 \beta w$ and $2z \equiv 2\alpha w$ (all mod 2). That is 
 $2{\bf u} \equiv 2w {\bf c} \bmod 2$, so ${\bf u} \in O$. (so far using only 
 (\ref{division}), (\ref{ab1}) and (\ref{maincongruence}))

Now let $O'$ be an order which contains ${\bf i}$ and ${\bf j}$, and suppose ${\bf v}=(x,y,z,w) \in O'$, then because $R$ is integrally closed,
${\rm tr}({\bf v})$,
${\rm tr}({\bf i}{\bf v})$, ${\rm tr}({\bf j}{\bf v})$, ${\rm tr}({\bf k}{\bf v})$ and ${\rm N}({\bf v})$ are all in $R$ (\cite{V}, Corollary 3.6). These give in turn $2x \in R$, $2ay \in R$, $2bz \in R$, 
   $2abw \in R$ and $x^2-ay^2-bz^2+abw^2 \in R$. Setting $X=2x$, $Y=2ay$, $Z=2bz$, $W=2abw$, the last equation gives 
\begin{equation}
\label{normwithclearedfractions}
abX^2-bY^2-aZ^2+W^2 \in 4abR
\end{equation}
whence
\begin{equation}
\label{getdivisibility}
W^2-bY^2 \in aR\;\;\;\;\;\;W^2-aZ^2 \in bR
\end{equation}
By (\ref{linind})
it follows that $a|Y$ and $b|Z$. Together with (\ref{normwithclearedfractions}) and (\ref{ab1}), this gives $ab|W$.
It follows that $x,y,z,w \in \frac{1}{2}R$, and so from the first statement that ${\bf v} \in O'$.\hfill  $\Box$

\bigskip

Theorem \ref{ax} with $R=\IZ$ and $a=b=-1$, $\alpha=\beta=1$ gives Theorem \ref{Hurwitzorder}. 
In this case (\ref{maincongruence}) is the statement that $\IZ$ has the four squres property.

\bigskip

\noindent{\bf Proof of Theorem \ref{halfintgroup}}. First we show that any element of ${\cal O}$ has half-integer coefficients. 
Let ${\bf u}=(R,S,T,W) \in {\cal O}$, let $d$ be the lowest common denominator of all the coefficients
(reduced as far as possible) of the components of ${\bf u}$, then $d{\bf u} \in \IZ[u,v]^4$ and $N(d{\bf u}) \in d^2 \IZ[u,v]$. If $d$ is divisible by an odd prime $p$, then reducing the coefficients in
$d{\bf u}$ mod $p$ we obtain a nonzero quaternion in ${\cal Q}_p:=\left(\frac{u^2-1,v^2-1}{\IZ_p[u,v]}\right)$,
which has zero norm, but this is impossible as we will show that ${\cal Q}_p$ is a division algebra.
By \cite{MR}, Theorem 2.3.1, it suffices to show that the equation
\[
(u^2-1)p^2(u,v)+(v^2-1)q^2(u,v)=1
\]
has no solution with $p,q \in \IZ_p(u,v)$. Setting $v=1$ this equation becomes $(u^2-1)p^2(u,1)=1$,
which clearly has no solution, as $(u^2-1)$ is not a square. So we conclude that $d$ is a power of 2.  If $d$ were a multiple of 4, then $d{\bf u}$ would have integer coefficients, not all even, and norm divisible by 16, but the second part of Corollary \ref{ax2congruences}, with $a=u^2-1$, $b=v^2-1$, 
$\alpha=u+1$, $\beta=v+1$, $\alpha'=u-1$, $\beta=v-1$, shows that
this is impossible.

To complete the proof, we apply Theorem \ref{ax} with
$R=\IZ[u,v]$, $a=u^2-1$, $b=v^2-1$, $\alpha=u+1$, $\beta=v+1$. In this case, we can easily verify 
(\ref{ic})-(\ref{ab1}).
To prove (\ref{linind}), let $u^2-1 | p^2(u,v)+(1-v^2)q^2(u,v)$, 
where $p,q \in \IZ[u,v]$. For all $v \in (-1,1)$, both summands on the right hand side are nonnegative. Hence, when $u=\pm 1$, both vanish. It follows that $u^2-1$ divides $p$ and $q$.
Together with the corresponding statement obtained by interchanging $u$ and $v$, this gives (\ref{linind}).
Finally, the first part of Corollary \ref{ax2congruences} gives (\ref{maincongruence}). 
\hfill $\Box$

\bigskip

\begin{lemma} There is a member of $u\in {\cal U}_1$ with irrational coefficients.
\label{irrational}
\end{lemma}
An example is the quartic
\begin{eqnarray*}
{\bf u}&=&[(1-u^2)(a-av^2+v^2)+ u^2v,  (v-1)((b-au)(v+1)+uv),\\
&& (1-a)(1-v)(1-u^2)+u,  a+bv-u(a-1)(v-1)],
\end{eqnarray*}
which has norm 1 whenever $2a - 3a^2 - b^2=1 - 2a +  a^2 - ab=0$. A routine calculation shows that these have real solutions $a=b=1/2$, and where $a$ and $b$ are the (unique) real roots of $2x^3 - 2x^2 + 2x - 1$ and $2x^3 +6x^2 + 4x - 1$ respectively.  These roots are not rational.

\bigskip

\subsection{Generation.}

Here we consider the question as to whether or not ${\cal U}$ finitely generated.  We thank Alan Reid for providing us with a simpler proof than our earlier argument based on arithmetic Kleinian groups.

\begin{theorem}  The group ${\cal U}$ is not finitely generated.  \end{theorem}
\noindent{\bf Proof.} 
We may identify ${\cal U}$ is the obvious way with the group of elements of norm 1 in the quaternion algebra
\begin{equation}\label{quatZ}
{\cal Q}_\IQ=\left( \frac{u^2-1, v^2-1}{\IQ(u,v)} \right),
\end{equation}

Suppose that $O$ is a (maximal) order in ${\cal Q}_\IQ$
and simply specialize $u,v$  as follows. Put $u=0$ and,  for $p$ a prime $p\geq 3$,
\[ v = \frac{1}{2} \big(p+\sqrt{p^2-4}\big), \]
which has conjugate $\bar v = \frac{1}{2} \Big(p-\sqrt{p^2-4} \Big)\in (-1,1)$.  Then ${\cal Q}_\IQ$ has homomorphic image
\begin{equation}
\left( \frac{-1,\frac{1}{2} (p+\sqrt{p^2-4} )}{\IQ(\sqrt{p^2-4}) }\right),
\end{equation}
Apart from the identity,  the other real embedding is $\sigma(v)=\bar v$ and so $\sigma(v^2-1)=\bar v^2-1 < 0$.  Hence
the group of elements of norm 1 in the order $O$ so specialized
is some arithmetic Fuchsian group coming from a division algebra over
${\cal Q}_\IQ(v)$,  see Theorem \ref{arith} below.  

Now the rank of this group must go to infinity with $p$ as there are only finitely many arithmetic Fuchsian groups whose quotients are surfaces of a given topological type,  \cite{MR}. In particular,  this implies the group of elements of norm 1 of $O$ cannot be finitely generated. \hfill $\Box$

\bigskip

Calculation shows that the 5-element set
\begin{eqnarray}
\label{gens}
\left\{(u, 1, 0, 0), \;(v, 0, 1, 0), \; (u, v, 0, -1),  \nonumber \frac{1}{2}(1 + u + v - u v, v-1, 1-u, -1),\right.\\  \left. \frac{1}{2}(1 + u + v - u v, v - 1, 1 - u, 1) \right \},
\end{eqnarray}
each of which is of degree 1 or 2, generates every member of ${\cal U}$ of degree at most 4. However we also have for example (proof omitted) that ${\bf u}$ below does not lie in the subgroup generated by these elements.
\begin{eqnarray*}
{\bf u}&=& \frac{1}{2}( -1 + u^2 - 2 u^3 - v^2 + 3 u^2 v^2 + 2 u^3 v^2,  1 - u + 2 u^2 - v^2 + 5 u v^2 - 2 u^2 v^2, \\
&&  1 - u^2 + v - 2 u v + u^2 v + 4 u^3 v,   1 + u - v + 3 u v - 4 u^2 v)
\end{eqnarray*}

\subsection{Quaternions as Isometries}
\label{asisometries}

We now look at what happens to the members of ${\cal U}$ under the evaluation map
$\psi_{\lambda, \lambda',\mu}$. We will assume for the moment that $\lambda, \mu \notin [-1,1]$ and $\lambda'=-1$, and abbreviate $\psi_{\lambda, -1,\mu}$ to $\psi$. We set $\Gamma=\Gamma(\lambda, \mu)=\psi({\cal U})$. Now (\ref{newad}) and (\ref{newoffdiagonals}) become
\begin{equation}
\label{newadspecial}
a=-d=\frac{1}{\sqrt{2}}\sqrt{-1-\mu}, \;\;\;\;c=-b= \frac{1}{\sqrt{2}} \sqrt{1-\mu}
\end{equation}
and in addition
\[
ab=cd=\frac{1}{2} \sqrt{1-\mu}\sqrt{-1-\mu}=\pm \frac{1}{2} \sqrt{\mu^2-1},
\]
and so $\psi((R,S,T,W))$ is
\begin{equation}
\label{matrixform3}
\left[\begin{array}{cc}
R_w+S_w \sqrt{\lambda^2-1} & \pm(T_w +W_w \sqrt{\lambda^2-1})\sqrt{\mu^2-1}  \\ \pm(T_w -W_w \sqrt{\lambda^2-1})\sqrt{\mu^2-1} & R_w -S_w \sqrt{\lambda^2-1}
\end{array}
\right].
\end{equation}

First we revisit the three quaternions ${\bf w}_i \in {\cal V}_0$ ($i=1,2,3$) defined at (\ref{ws}), which have images $\tilde{{\bf w}}_i:=\rho({\bf w_i})\in {\cal V}$, namely
\begin{eqnarray}
\tilde{{\bf w}}_1&=& (u,1,0,0), \nonumber \\
\tilde{{\bf w}}_2&=& (u,v,0,-1), \nonumber\\
\label{newws}
\tilde{{\bf w}}_3&=& \frac{1}{2}(1+u+v-uv, v-1, 1-u,-1),
\end{eqnarray}

As we have already seen (or directly from (\ref{matrixform3})), these map respectively to
the isometries $f^2$, $gf^2g^{-1}$ and $[g,f]$, where $\lambda(f)=\lambda$, $\mu(f,g)=\mu$, $f$ is loxodromic (since $\lambda \notin [-1,1]$), $\ax(f)=(0,\infty)$, and $g$ is an order 2 elliptic whose axis is disjoint from $\ax(f)$ (since $\mu \notin [-1,1]$), and has mutually reciprocal endpoints.
As noted in the remarks after the proof of Theorem \ref{canonicalgroup}, this means that the common perpendicular of  $\ax(f)$ and $\ax(gf^2g^{-1})$ has endpoints $\pm 1$.

\bigskip

We have now got back the subgroup of $\Gamma$ comprising the
balanced even words in $f$ and $g$ (since $g$ is order 2, the distinction between regular and irregular words now vanishes). We can now extend this subgroup. Let $\varphi_f (z)=-z$, $\varphi_h (z)=1/z$ $\varphi (z)=-1/z$; these three isometries are each of order 2, have mutually orthogonal axes and generate a Klein 4-group, $K$. We define
$h=g \varphi_f$. Recall (\ref{newAB}) that $g$ has matrix representative $B=\left[\begin{array}{cc}
a & b \\ c & d
\end{array}
\right]$. Thus, using (\ref{newadspecial}), $h$ and $h^2$ have respective matrix representatives
\[
M_h=\frac{1}{\sqrt{2}}\left[ \begin{array}{cc} \sqrt{\mu+1} & \pm \sqrt{\mu-1} \\ \pm \sqrt{\mu-1} & \sqrt{\mu+1} \end{array} \right], \;\;\;\;\;
M_{h^2}=\left[ \begin{array}{cc} \mu & \pm \sqrt{\mu^2-1} \\ \pm \sqrt{\mu^2-1} & \mu \end{array} \right]
\]
 The axis of $h$ has endpoints $\pm 1$, and $\lambda(h)=\mu$. Also $h^2 \in \Gamma$; specifically $h^2=\psi((v,0,\pm 1,0))$.

At this point a certain symmetry between $f$ and $h$ is becoming apparent. Both are loxodromic, both have squares in $\Gamma$ and their axes are mutually perpendicular. To develop this symmetry further we express $h$, like $f$, as a product of two order 2 elliptics.
Set $\tilde{g}=f \varphi_h$; explicitly, $\tilde{g}(z)=A/z$, where $f(z)=Az$, so that $\tilde{g}$ has order 2. We now have 
\begin{eqnarray*}
g \varphi_f=h, && \tilde{g} \varphi_h=f,\\
\lambda(h)=\mu(f,g), && \lambda(f)=\mu(h,\tilde{g}),\\
\ax(\varphi_f)=\ax(f), && \ax(\varphi_h)=\ax(h).
\end{eqnarray*}
We can summarise all this by saying that the pair $(h,\tilde{g})$ is obtained from $(f,g)$ (up to conjugacy) by interchanging the parameters $\lambda$ and $\mu$.
\begin{theorem}
\label{perpaxes}
The subgroup $P$ of $\langle f,h \rangle$ comprising the isometries of the form $f^{n_1}h^{m_1}f^{n_2}h^{m_2}\ldots f^{n_k}h^{m_k}$, where $n_1+n_2+ \ldots n_k$ and $m_1+m_2+ \ldots m_k$ are both even,
is a subgroup of $\Gamma$.
\end{theorem}

\noindent{\bf Sketch of Proof.} We first show that $P=\langle f^2, h^2, fh^2f^{-1}, hf^2h^{-1} \rangle$. This can be done using induction along the same lines as the proof of Lemma \ref{evengen}. We have already seen that $f^2, h^2 \in \Gamma$.
The proof is completed by showing that $fh^2f^{-1}$, $hf^2h^{-1}$ have respective matrix representatives
\begin{eqnarray*}
M_{fh^2f^{-1}}& = & \left[ \begin{array}{cc} \mu & \pm (\lambda+\sqrt{\lambda^2-1})\sqrt{\mu^2-1} \\ \pm (\lambda-\sqrt{\lambda^2-1})\sqrt{\mu^2-1} & \mu \end{array} \right]\\
& = & \psi((v,0, \pm u ,\pm 1))\\ \\
M_{hf^2h^{-1}}& = & \left[ \begin{array}{cc} \lambda+\mu \sqrt{\lambda^2-1} & \mp \sqrt{\lambda^2-1}\sqrt{\mu^2-1} \\
\pm \sqrt{\lambda^2-1}\sqrt{\mu^2-1} & \lambda-\mu \sqrt{\lambda^2-1} \end{array} \right] \\ & = & \psi((u,v, 0 ,\mp 1)).
\end{eqnarray*}
\hfill $\Box$

Clearly $P$ is a finite-index subgroup of $\langle f,h \rangle$, and it follows in particular that if $\Gamma$ is discrete, then so is $\langle f,h \rangle$.
Further (see (\ref{matrixform3})) a sufficient condition for this is that $\lambda$ and $\mu$ both lie in a discrete subring of $\IC$ (i.e. a subring of the ring of integers of some imaginary quadratic field).

\begin{corollary}
\label{perpaxes}
If $R$ is a discrete subring of $\IC$, $f$ and $h$ are non-parabolic non-identity isometries in ${\rm Isom}^+(\IH)$ with perpendicular axes, and $\lambda (f), \lambda (h) \in R$, then $\langle f,h \rangle$ is discrete.
\end{corollary}
In particular we have discreteness when $\lambda (f), \lambda (h)$ are integers. Another discrete example is
$\lambda (f)=\lambda (h)=\frac{-1+\sqrt{3}}{2}$, which minimizes  $\max\{\tau_f,\tau_h\}$ among all 
two generator non-elementary groups having loxodromic generators with perpendicular axes,   \cite{MM2}. We discuss discreteness criteria further in Section \ref{arithmetic}. 

\medskip

Additionally,  we can add in all the order 2 elliptics and preserve discreteness. These elliptics fall into three Klein 4-groups: $K:=\{\varphi_f, \varphi_h, \varphi\}$, $K_f:=\{\varphi_f, \tilde{g}, \tilde{g} \varphi_f  \}$ and $K_h:=\{\varphi_h, g, g \varphi_h \}$.
\begin{theorem}
\label{addelliptics}
The group $P_1$ generated by $K$ $K_f$ and $K_h$ is an extension of
$\langle f,h \rangle$ of index at most 2.
\end{theorem}
{\bf Proof.} Every $\alpha \in P_1$ can be represented by a word in $\{\varphi_f, \varphi_h, \varphi, g, \tilde{g},f,h\}$
which we suppose to have the fewest possible elliptic letters, and with the first elliptic letter occurring as close to the right as possible.
If $a \in K \cup K_f$ then
$afa$=$f^{\pm 1}$, so that $af=f^{\pm 1}a$, and similarly $af^{-1}=f^{\mp 1}a$. If $a=g$ then $af^{\pm 1}=gf^{\pm 1}=g(\varphi_f f^{\pm 1} \varphi_f^{-1})g^{-1}g=hf^{\pm 1}h^{-1}g$, because $f$ and $\varphi_f$ commute. 
By our assumptions about the word, it follows that no elliptic letter can immediately precede an $f^{\pm 1}$, and symmetrically it cannot immediately precede an $h^{\pm 1}$ either. It follows that all of the elliptic letters are at the right of the word. But the product of any two elliptics is either another elliptic or a product of (at most two) of the loxodromics $f^{\pm 1}$ and $h^{\pm 1}$, so the word contains at most one elliptic letter. Thus $P_1$ is an extension of $\langle f,h \rangle$ of index at most 2.
\hfill $\Box$

\section{Arithmeticity.} \label{arithmetic}

In this section we first recall some further terminology concerning quaternion algebras with an aim to extending the discreteness conditions described above.  This section is adapted from \S 4 of \cite{GMMR}.

\bigskip

Let $k$ be a number field. A {\em place} $\nu$ of $k$ is an equivalence class of valuations on $k$. Such a place is real (complex) if it is associated to a real embedding (conjugate pair of complex embeddings) of $k$. We denote by $k_\nu$ the completion of $k$ at the place $\nu$. If $Q$ is a quaternion algebra over $k$, we say that $Q$ is ramified at $\nu$ if $Q\otimes_k k_\nu$ is a division algebra of quaternions. Otherwise $\nu$ is unramified. If $\nu$ is a real place, then $Q$ is ramified if and only if $Q\otimes_k k_\nu\equiv \IHM$. It is straightforward to check whether a quaternion algebra $Q = \left( \frac{a,b}{k} \right)$ is ramified at a real place $\nu$; if $\nu$ corresponds to
the real embedding $\sigma$, then $Q$ is ramified at $\nu$  if and only if $\sigma(a)$ and $\sigma(b)$ are both negative.

\medskip
We can now define an arithmetic Kleinian group. Let $k$ be a number field with one complex place and $Q$ a quaternion algebra over $k$ ramified at all real places. Next let $\rho$  be an embedding of $Q$ into $SL(2,\IC)$, let $O$ be an  order of $Q$ and $O_1$ the elements of norm $1$ in $O$. Then $\rho(O_1)$ is a discrete subgroup of $SL(2,\IC)$ and its projection to $PSL(2,\IC)$  is an  Kleinian group.  Kleinian groups so constructed, together with those which are commensurable to them, are {\em arithmetic}.
  We note in passing that arithmetic Fuchsian groups arise in a similar manner. However in that case, the number field is totally real and the algebra ramified at all real places except the identity.

\medskip
For a subgroup $\Gamma$ of $SL(2,\IC)$ the invariant trace field is defined as
\begin{equation}
k\Gamma = \IQ(\{\tr^2(g):g\in \Gamma\})
\end{equation}
Then we set
 \[
 Q\Gamma=\left\{
\sum a_i\;g_i : a_i \in \IQ(\tr(\Gamma)), g_i \in \Gamma  \right\} \]
Then $Q\Gamma$ is a quaternion algebra over $\IQ(\tr(\Gamma))$.  Additionally, if $ \tr(\Gamma)$ consists of algebraic integers we see that
 \[ O\Gamma=\left\{ \sum a_i\;g_i : a_i \in R_{\IQ(\tr(\gamma))}, g_i \in \Gamma \right\} \]
is an order in $Q\Gamma$.  Here  $R_{\IQ(\tr(\gamma))}$ is the ring of integers in $\IQ(\tr(\Gamma))$.
Then $\Gamma$ is arithmetic if and only if the following conditions are satisfied:
 \begin{enumerate}
\item  $k\Gamma$ is an algebraic number field;
\item  $\tr(\Gamma)$ consists of algebraic integers;
\item for every $\IQ$-isomorphism $\sigma : k\Gamma \to  \IC$, other than the identity or complex conjugation,  $\sigma(\tr(\Gamma^{(2)}))$ is
bounded in $\IC$.
 \end{enumerate}

 \medskip

In practice, it is hard to apply this characterisation directly, the problem being to establish the boundedness of the traces at real embeddings. However, in \cite{GMMR} we obtained the following more useful method for proving groups discrete.
\begin{theorem}\label{arith} Let $\Gamma$ be a finitely generated non-elementary subgroup of $SL(2,\IC)$ such that
  \begin{enumerate}
\item  $k\Gamma$ has exactly one complex place or is totally real;
\item $\tr(\Gamma)$ consists of algebraic integers;
\item $Q\Gamma^{(2)}$ is ramified at all non-identity real places of $k\Gamma$,
\end{enumerate}
then $\Gamma$ is a
subgroup of an arithmetic Kleinian or Fuchsian group.
\end{theorem}

\begin{corollary}  A group of elements ${\cal G}$ of norm 1 in an order ${\cal O}$ of the quaternion algebra ${\cal Q}_\IQ$,
\[ {\cal Q}_\IQ=\left( \frac{u^2-1, v^2-1}{\IQ(u,v)} \right) \]
is a discrete subgroup of an arithmetic Kleinian group if
\begin{itemize}
\item $u$ is a complex root of a monic irreducible polynomial of degree $n$ with integer coefficients $p(z)$ all of whose $(n-2)$ real roots $r_1,r_2,\ldots,r_{n-2}$ lie in the interval $(-1,1)$.
\item $v$ is an algebraic integer in $\IQ(u)$.
\item for each non-identity real embedding $\sigma_i:\IQ(u)\to\IQ$, $\sigma|\IQ=id$,  defined by $\sigma_i(u)=r_i$,   the image $\sigma(v)\in (-1,1)$.
\end{itemize}
\end{corollary}
\noindent{\bf Proof.}  The first condition gives $\IQ(u)$ a number field of degree $n$ over $\IQ$ with one complex place, and $v\in \IQ(u)$ then gives $\IQ(u,v)=\IQ(u)$.  If $\sigma_i$ is a real embedding, then  $\sigma_i(u^2-1)=\sigma_i(u)^2-1=r_{i}^2-1<0$ and by hypothesis $\sigma_i(v^2-1)<0$ so the quaternion algebra is ramified at all the real places.   Next,  the trace is $2R(u,v)\in \IZ(u,v)$. This must be an algebraic integer under our hypotheses. \hfill $\Box$

\section{Discreteness: necessary conditions.}
Let  $\Gamma:=\psi_{\lambda, \lambda',\mu}({\cal U})$, a subgroup of $SL(2, \IC)$ which, we recall, comprises the matrices of the form
\[
W=\left[\begin{array}{cc}
R(\lambda, \mu)+S(\lambda, \mu)\sqrt{\lambda^2-1} & 2ab(T(\lambda, \mu)+W(\lambda, \mu)\sqrt{\lambda^2-1})  \\ 2cd(T(\lambda, \mu)-W(\lambda, \mu)\sqrt{\lambda^2-1} )& R(\lambda, \mu)-S(\lambda, \mu)\sqrt{\lambda^2-1}
\end{array}
\right] 
\]
where the polynomials $R$, $S$, $T$ and  $W$ are polynomials with half-integer coefficients satisfying (\ref{symmdet}), with 
\[ a=\frac{1}{\sqrt{2}}\left(\sqrt{\lambda'+1}+\sqrt{\lambda'-\mu}\right),\;\; d=\frac{1}{\sqrt{2}}\left(\sqrt{\lambda'+1}-\sqrt{\lambda'-\mu}\right),\]
and $b$ and $c$ are given by (\ref{offdiagonals}) when $\mu=1$,
and $c=-b= \sqrt{\frac{1-\mu}{2}}$ otherwise.

We know that $\Gamma$
extends the group of regular even balanced words in $A$ and $B$ given by
\begin{equation}\label{AB3}  A= \left(
\begin{array}{cc}
 \frac{\lambda+\sqrt{\lambda^2-1}-1}{\sqrt{2} \sqrt{\lambda-1}} & 0 \\
 0 & \frac{-\lambda+\sqrt{\lambda^2-1}+1}{\sqrt{2} \sqrt{\lambda-1}} \\
\end{array}
\right), \;\;\;\;\; B=\left(
\begin{array}{cc}
a&b \\
c & d\\
\end{array}
\right), \end{equation}
 
 As we have observed (Theorem \ref{groupindex}) this group is an index two subgroup of the group of all regular balanced words in $A$ and $B$. Similarly we can show that $\Gamma$ an index two subgroup of the group $\tilde{\Gamma}$ generated by $\Gamma$ and $A$. To see this, it is enough to show that, for $G \in \Gamma$, $AGA^{-1} \in \Gamma$, a straightforward calculation using
(\ref{generalconjuation}). It follows in particular that $\Gamma$ is discrete if and only if $\tilde{\Gamma}$ is.

Another routine calculation, using (\ref{generaltrace}) and the facts that $\beta(A) =2(\lambda-1)$
and
\[4abcd = bc(1+bc) = 4 \Big(\frac{\mu-1}{2}\Big)\Big(1+\frac{\mu-1}{2}\Big) =  \mu^2-1, \]
gives
\begin{eqnarray*}
\gamma(A,W)& = & \tr[A,W]-2= -8(\lambda-1)abcd \big(T^{2}(\lambda, \mu)-W^2(\lambda, \mu)(\lambda^2-1)\big) \\
& = &  -2(\lambda-1)(\mu^2-1)\big(T^{2}(\lambda, \mu)-W^2(\lambda, \mu)(\lambda^2-1)\big)
\end{eqnarray*}
Since $R^2-(\lambda^2-1)S^2-(\mu^2-1)( T^2 -(\lambda^2-1)W^2)=1$ we then have
\[ \gamma(A,W) -\beta(A) =-2(\lambda-1)\big[R^2(\lambda,\mu)-(\lambda^2-1)S^2(\lambda,\mu) \big] \]
We write $\beta=\beta(A)$ and $\tilde{\gamma}=\gamma(A,W)$ to obtain the following three identities:
\begin{eqnarray*}
\label{id1}
|\beta|+|\tilde{\gamma}| & = &   2|\lambda-1|\Big(1+|\mu^2-1|\big|T^{2}(\lambda, \mu)-W^2(\lambda, \mu)(\lambda^2-1)\big| \Big) \;\;\;\;\;\;\; \;\;\\
 \label{id2} |\beta|+|\beta-\tilde{\gamma}| & = &  2|\lambda-1|\Big(1+ \big|R^2(\lambda,\mu)-(\lambda^2-1)S^2(\lambda,\mu)\big| \Big)
 \end{eqnarray*}
 and $|\tilde{\gamma}||\beta -\tilde{\gamma}|$
\[ = 4|\lambda-1|^2 |\mu^2-1|  \big|R^2(\lambda,\mu)-(\lambda^2-1)S^2(\lambda,\mu)\big|  \big|T^{2}(\lambda, \mu)-W^2(\lambda, \mu)(\lambda^2-1)\big| \]
The three equations above enable us to use the following test for
the discreteness of $\tilde{\Gamma}$ (and so of $\Gamma$).

\begin{theorem}
Let $c_0 = 2-2\cos(\pi/7) \approx 0.198062$ and $\lambda^2,\mu^2\neq 1$.  Then with the notation above, $\Gamma$  is discrete if and only if for every $W \in \Gamma$ the three inequalities
\begin{eqnarray}
|\beta|+|\tilde{\gamma}| & \geq  &  1, \hskip25pt \mbox{if $\tilde{\gamma} \neq 0$,  and} \\
|\beta|+|\beta -\tilde{\gamma}| & \geq  & 1, \hskip20pt \mbox{ if $\tilde{\gamma}\neq \beta$,  and} \\
|\tilde{\gamma}||\tilde{\gamma}-\beta| & \geq  & c_0, \hskip18pt \mbox{ if $\tilde{\gamma} \neq 0,\beta$.}
\end{eqnarray}
\end{theorem}
\noindent{\bf Proof.}  $[\Rightarrow]$ First suppose $\{W_i\}_{i=1}^{\infty} \subset {\cal U}$ is an infinite sequence,  that  $W_i\to id$ as $i\to \infty$,  and that (with the obvious notation) $\gamma_{i} \neq 0,\beta$.  Then of course ultimately the last inequality is violated since $\gamma_i=\gamma(A,W_i)=\tr[A,W_i]-2 \to 0$.  To remove the assumption that $\gamma_{i} \neq 0,\beta$ we consider two cases.

{\bf Case 1.}  $\gamma_i=0$ for infinitely many $i$.  Then $W_i$ has a fixed point in common with $A$ in $\oC$.  Now  $X=BAB^{-1} \in \Gamma$ (only if $B$ is order 2).  If $X$  shares a fixed point with $A$,  or maps one fixed point to another,  then $A$ and $XAX^{-1}$ have a common fixed point in $\oC$ and hence
\[ 0=\gamma(A,X)= \gamma(A,BAB^{-1}) =\gamma(A,B)(\gamma(A,B)-\beta) \]
so $\gamma(A,B)=0$ or $\gamma(A,B)=\beta$.  However,  $v=1-2\gamma/\beta\in \{\pm 1\}$ in either case,  and this is excluded by hypothesis.   We now deduce that $X^{-1}W_iX$ does not share a fixed point with $A$ for infinitely many $i$,  and $XW_iX^{-1} \to id$.

\noindent {\bf Case 2.} $\gamma_i=\beta$  for infinitely many $i$. Then  $V_i=W_i A W_i^{-1}\to A$, $0=\gamma(A,V_i)=\gamma(A,A^{-1}V_i)$ and so we reduce to the first case by replacing $W_i$ by $[A,W_i]\to id$.

\medskip

$[\Leftarrow]$ Next,  suppose the group ${\cal U}$ is discrete,  but one of these inequalities is violated for some $W \in {\cal U}$.  The first two inequalities are J\o rgensen's inequality and a well known variant of it \cite{GMit}.  These are necessary conditions for the discreteness of the group $\langle A,B\rangle$ provided this group is not virtually abelian.  The last condition is a result of Cao \cite{Cao} improving other versions of inequalities J\o rgensen found \cite{Jorg,GMHalmos} for discrete groups generated by two elements of the same trace.  We state this in the following lemma.

\begin{lemma}
If $\langle f,g\rangle$ is Kleinian and $\beta(f)=\beta(g)$,   then $|\gamma(f,g)|\geq c_0$ where
\[ c_0 = 2-2\cos(\frac{\pi}{7}) \]
This bound is sharp and achieved in the $(2,3,7)$-triangle group.
\end{lemma}

Thus the violation of one of these inequalities shows that $\langle A,W\rangle$ is virtually abelian.  If the group is abelian,  then $\gamma(A,W)=0$.  If $WAW^{-1}=A^{-1}$,  the dihedral case,  then
$WAW^{-1}A^{-1}=A^{-2},$ and hence
\[ \gamma(A,W) =  \tr(WAW^{-1}A^{-1})-2=\tr A^2 -2 =\tr^2(A)-4 =\beta(A)\]
If $A$ is loxodromic,  then $\langle A,W \rangle$,  if discrete, is Kleinian unless $W$ fixes or interchanges the fixed points of $A$.  Otherwise there would be three,  and hence uncountably many limit points,  \cite{B}.  These reduce to the cyclic or dihedral cases.   If $A$ is elliptic,  then $|\beta|<1$ is required to violate either of the first two inequalities.  That is $A$ has order $7$ or more.  The classification of the elementary discrete groups \cite{B} shows this to reduce to the abelian or dihedral cases as well.  What remains is the case $\langle A,W\rangle$ is a discrete group with the last inequality violated.  Then this group is elementary and as $\gamma(A,WAW^{-1}) = \gamma_w(\gamma_w-\beta)$ a little argument using the classification of the elementary discrete groups reduces to the previous cases.  \hfill $\Box$

\section{Examples}
\label{Examples}

 We calculate some of the polynomials for balanced, even, good words in $f$ and $g$, and investigate when these have the same axis as $f$.
For $W = fgf^5g^{-1}fgf^2g^{-1}f^{-3}$,
\begin{eqnarray*}
2r(u,v)&=&-1 + 3 u - 2 u^2 - 10 u^3 + 4 u^4 + 8 u^5 - v - 3 u v +
   8 u^2 v + 4 u^3 v \\ &&  - 8 u^4 v + 2 v^2  - 6 u v^2 - 6 u^2 v^2 +
   14 u^3 v^2 + 4 u^4 v^2 - 8 u^5 v^2 \\
2s(u,v)&=&  -1 + 2 u + 10 u^2 - 4 u^3 - 8 u^4 - v - 4 u v + 4 u^2 v + 8 u^3 v +
 2 u v^2 \\ &&  - 6 u^2 v^2  - 4 u^3 v^2 + 8 u^4 v^2 \\
2t(u,v)&=& (u-1) (-1 + 2 u + 4 u^2) (-1 - 2 u + 4 u^2 + 4 u^3 - 4 u v +
   4 u^3 v)\\
2w(u,v)&=& (1 + 6 u - 4 u^2 - 20 u^3 + 8 u^4 + 16 u^5 - 2 v + 6 u v + 8 u^2 v -
  20 u^3 v \\ && - 8 u^4 v + 16 u^5 v)
\end{eqnarray*}
The only solution of $t(u,v)=w(u,v)=0$ is $u=-1/2$, $v=-1/3$. For these values we also have $s(u,v)=0$ and $r(u,v)=-1$, i.e. $W$ is a relator of $\langle f,g \rangle$ for these values.  However,  we return back to (\ref{lambda}) and (\ref{mu}) to see
\[ \beta = 2u-2 = -3,  \;\;\;\;\; \gamma = \beta(1-v)/2= -2 \]
so $f$ has order three and $f$ and $g$ .
\begin{corollary}  Let $\Gamma$ be a Kleinian group and $f,g\in \Gamma$.  Then
\[  fgf^5g^{-1}fgf^2g^{-1}f^{-3}=1 \]
 if and only if $f$ has order $3$ and $\langle f,g \rangle$ is a Euclidean triangle group or an abelian group.
\end{corollary}
Of course if $f$ and $g$ commute then $ fgf^5g^{-1}fgf^2g^{-1}f^{-3}=1$,  then $ f^6  =1$.

\bigskip

For $W = fgf^5 g^{-1}f^{-2}$
\begin{eqnarray*}
2r(u,v)&=&-1 - 3 u + 2 u^2 + 4 u^3 - v + 3 u v + 2 u^2 v - 4 u^3 v \\
2s(u,v)&=&  1 + 2 u - 4 u^2 - v + 2 u v + 4 u^2 v \\
2t(u,v)&=& (u-1) (1 + 2 u) (-1 + 2 u + 4 u^2) \\
2w(u,v)&=& (-1 + 2 u) (-1 + 2 u + 4 u^2)
\end{eqnarray*}
This time $t$ and $w$ have a common factor $-1 + 2 u + 4 u^2$, so that they vanish simultaneously when $u = 1/4 (-1 \pm \sqrt{5})$, and for all values of $v$. However
$s(1/4 (-1 \pm \sqrt{5}),v)=1/2 (-1 + \sqrt{5}) \neq 0$, so $W$ can never be a relator of $\langle f,g \rangle$.
 \bigskip

For $W = fgfg^{-1}f^2gfg^{-1}fgf^{-1}g^{-1}f^2gf^{-1}g^{-1}$
\begin{eqnarray*}
r(u,v)&=&-u - u^2 + u^3 + u^4 + u^5 - v^2 + u v^2 - u^2 v^2 + u^3 v^2 +
 2 u^4 v^2 \\
 &-& 2 u^5 v^2 + v^4 - 3 u v^4 + 2 u^2 v^4 + 2 u^3 v^4 -
 3 u^4 v^4 + u^5 v^4 \\
s(u,v)&=&  -1 + u^2 + 2 u^3 + u^4 + v^2 - 2 u v^2 + 3 u^2 v^2 -
 2 u^4 v^2 - v^4 \\
 &+& 2 u v^4 - 2 u^3 v^4 + u^4 v^4 \\
t(u,v)&=& (u-1)(u+1) (u - v + u v) (-1 - u - u^2 - v + u v + v^2\\
 &-& 2 u v^2 + u^2 v^2) \\
w(u,v)&=& (u-1) (-u - u^2 - u^3 - v - u v - u^2 v - u^3 v + v^2 - u v^2 -
   u^2 v^2 \\
  &+& u^3 v^2 + v^3 - u v^3 - u^2 v^3 + u^3 v^3)
\end{eqnarray*}
Here $W$ is a relator in the group that minimizes the maximum of the two translation lengths $\max\{\tau_f, \tau_h\}$, when $f$ and $h$ are two loxodromics with perpendicular axes \cite{MM2}.
Now $t$ and $w$ have a common factor $(u-1)$, so that they vanish simultaneously when $u = 1$. However
$t=w=0$ also holds when $v=0$ (perpendicular axes), when $u=0$ and when $u=-1\pm i\sqrt{3}$. In the last two cases $r=1$ and (consequently) $s=0$.

These examples raise some general questions:
\begin{enumerate}
\item Which words can be relators? (ie for which words do $s=t=w=0$, $r=\pm1$ have a solution, apart from the trivial solutions $u=1$ ($f=Identity$) and $v=1$ ($f$ and $g$ have the same axis)?)
\item For which words do $t$ and $w$ have a non-constant common factor in which neither of the variables $u$ and $v$ is absent? (giving an infinite family of solutions for $t=w=0$)
\item For which words do $t$ and $w$ have no such common factor, so that $t=w=0$ has only finitely many roots, and one of these roots also makes $s=0$ and (hence) $r=\pm 1$ (ignoring the trivial cases $u=1$, $v=1$). Is any such group discrete?
\end{enumerate}

\subsection{Explicit Formulae}
\label{explicit}
We now give (without going into details of the computation) explicit values for the polynomials $R=R_w$, $S=S_w$, $T=T_w$ and $W_w$ of (\ref{matrixform2}) associated with the even word $w=f^{n_1}gf^{n_2}g^{-1}f^{n_3}gf^{n_4}g^{-1}f^{n_5}$, which are expressed in terms of Chebyshev polynomials indexed by various combinations of the powers $n_i$. For such a word, and for $S \subseteq \{1,2,3,4,5\}$, we let
$T_S(u)=T_{(\epsilon_1 n_1+\epsilon_2 n_2 +\epsilon_3 n_3+\epsilon_4 n_4+\epsilon_5 n_5)/2}(u)$,
$U_S(u)=U_{(\epsilon_1 n_1+\epsilon_2 n_2 +\epsilon_3 n_3+\epsilon_4 n_4+\epsilon_5 n_5)/2-1}(u)$,
where $\epsilon_i=-1$ if $i \in S$, $\epsilon_i=1$ if $i \notin S$,
e.g. $T_{\{2,3\}} (u)=T_{(n_1-n_2-n_3+n_4+n_5)/2}(u)$. (We set $U_{-1}(x)=0$, and, for $n<0$, $T_n(x)=T_{|n|}(x)$, $U_{n-1}(x)=-U_{|n|-1}(x)$). 

We have calculated:
\begin{eqnarray}
\label{integerindices}
R(u,v)&=& \frac{1}{4} \left( \sum_{S \subseteq \{2,3,4\}} {(-1)}^{|S|} T_S(u) \right)v^2 \nonumber \\
&+& \frac{1}{2} \left( T_{\emptyset}(u)-T_{\{2,4\}} (u)\right)v \nonumber \\
&+& \frac{1}{4}  \left( T_{\emptyset} (u)+T_{\{2\}}(u)+T_{\{3\}} (u)+T_{\{4\}}(u)+T_{\{2,4\}} (u)+T_{\{2,3,4\}} (u) \right. \nonumber 
\end{eqnarray}
$S(u,v)$ is the same, but with $U_S$ substituted for $T_S$ throughout in (\ref{integerindices}) above,
\begin{eqnarray*}
T(u,v)&=& \frac{1}{4} \left( \sum_{S \subseteq \{2,3,4\}} {(-1)}^{|S|} T_{S \cup\{5\}}(u) \right)v \\
&+& \frac{1}{4} \left( T_{\{3,5\}} (u)-T_{\{1,3\}} (u)+T_{\{2,5 \}} (u)-T_{\{1,4\}} (u) \right.\\
&-&  \left. T_{\{1,2,3\}} (u)+T_{\{3,4,5\}} (u)-T_{\{1\}} (u)+T_{\{5\}} (u)
\right).  
\end{eqnarray*}
\begin{eqnarray*}
W(u,v)&=& \frac{1}{4} \left( \sum_{S \subseteq \{2,3,4\}} {(-1)}^{|S|} U_{S \cup\{5\}}(u) \right)v \\
&+& \frac{1}{4} \left( U_{\{1,3\}} (u)+U_{\{3,5\}} (u)+U_{\{1,4 \}} (u)+U_{\{2,5\}} (u) \right.\\
&+&  \left. U_{\{1,2,3\}} (u)+U_{\{3,4,5\}} (u)+U_{\{1\}} (u)+U_{\{5\}} (u) \right).
\end{eqnarray*}

{\bf Remark:} These formulae exhibit the general fact that if the sequence $(n_1,n_2, n_3, n_4, n_5)$ is reversed, the sign of $t$ is changed, and $r$, $s$, $w$ are unchanged.

We thus have, for the shorter word $A^{n_1}BA^{n_2}BA^{n_3}$ ($n_4=n_5=0$)
\begin{eqnarray*}
R(u,v)&=& \frac{1}{2} \left[ \left( T_{(n_1+n_2+n_3)/2}(u)-T_{(n_1-n_2+n_3)/2}(u)\right)v \right.\\
&+&\left. \left( T_{(n_1+n_2+n_3)/2}(u)+T_{(n_1-n_2+n_3)/2}(u)\right)\right],
\end{eqnarray*}
$S(u,v)$ the same, but with $U_{n-1}$ substituted for $T_n$ throughout, and
\begin{eqnarray}
\label{linearinv}
T(u,v)&=& \frac{1}{2} \left[ T_{(n_1+n_2-n_3)/2}(u)-T_{(n_1-n_2-n_3)/2}(u) \right],
\end{eqnarray}
$W(u,v)$ the same, but with $U_{n-1}$ substituted for $T_n$ throughout.

\section{Roots of Trace polynomials.}

The purpose of this section is to establish a theorem which shows that the zero sets of the ``good word'' trace polynomials discussed in \S\ref{goodwords} are dense in the complement of the space of discrete and faithful representations of $\IZ_p*\IZ_2$ for $3\leq p \leq \infty$.  Indeed we show that the complement of the representations which are discrete and free on marked generators is the filled in Julia set of the semigroup of good word polynomials.

\begin{theorem} Let $f,g$ be M\"obius transformations with $\beta=\beta(f)\neq-4$, $\beta(g)=-4$, $\gamma=\gamma(f,g)$ and suppose that $\langle f,g\rangle$ is not discrete and free on the two generators $f$ and $g$.  Then for any open set $U$,  $\gamma\in U\in\IC$ there is a good word $w=w(f,g)$ for which the polynomial $q_w(z)=p_w(z,\beta)$, given by Theorem \ref{tracepolynomials}, has a root in $U$.
\end{theorem}
\noindent{\bf Proof.} There are two cases.

\subsection{$\langle f,g\rangle$ is discrete but not free on generators.} In this case there is a nontrivial word $w \in \langle f,g\rangle$ representing the identity.   Then,  for suitable $a$ one of $w,wf^a,f^aw$ is a good word as $g$ has order two.  Suppose this word is $v$.  Then $w=identity$ gives us  $0=\gamma(w,f)=\gamma(wf^a,f)=\gamma(f^aw,f) = p_v(\gamma,\beta)$ and so $\gamma$ itself is the root of a good word polynomial.  

\subsection{$\langle f,g\rangle$ is not discrete.} Let $U$ be a neighbourhood of $\gamma$ and define the good word polynomial zero set as  
\[ {\cal Z}=\{ z\in \IC\; : \; \mbox{ there is a good word $w$ so that $p_w(z,\beta)=0$ } \} \]
In \S \ref{examples} we gave a few examples of good words.  From that table we quickly deduce that among many other points
\[\{0,\beta, 1+\beta, 2+\beta \} \subset{\cal Z} . \]
Thus ${\cal Z}$ contains at least three finite points.  To simplify notation we suppress the $\beta$ variable in our word polynomials.  Next,  suppose that  for some good word $v$ we have $p_v(U)\cap {\cal Z}\neq \emptyset$.  Then there is a word $w\in \langle f,g\rangle$ and $z\in U$ such that $p_v(p_w(z))=0$.  However we know that the set of good words is closed under composition,  and $p_v(p_w(z))=p_{v*w}(z)$.  This point $z$ establishes the result.  We are left to consider the subcase 

\subsubsection{For all good words $v$, $p_v(U)\cap {\cal Z}= \emptyset$.}

So the family of analytic functions  $F_U=\{p_v|U:\mbox{ $v$ is a good word }\}$ has the property that  each element omits ${\cal Z}$ which contains at least three points.  Thus Montel's criterion shows that $F_U$ is a normal family.  If $\langle f,g\rangle$ is not discrete,  then there is a sequence of words $\{w_i\}_{i=1}^{\infty}$ in $\langle f,g\rangle$ with $w_i\to identity$ as $i\to\infty$ (this convergence is in the topology of $PSL(2,\IC)$,  that is in each entry of representative matrices).  Again,  as $g$ has order two,  $w_i$ are good words and we must have for any $z_0\in U$
\[ p_{w_i}(z_0) = \tr [f,w_i]-2 \to  \tr [f,identity]-2 = 0 \]
It follows that $p_{w_i}\to 0$ locally uniformly in $U$.  With our earlier argument using the fact that the set of good word polynomials is closed under composition,  we will be done if we can establish the density of the roots in some small neighbourhood of $0$. 

\subsubsection{Density of roots near $0$.}

We analyse this case in a fairly general framework using some of the theory of the dynamics of polynomial semigroups.  Much more can be found about this subject,  see for instance \cite{HM1,HM2,SS, sumi1,sumi2} and the references therein.  The point here is independent of $\beta=\beta(f)$,  we can find a good word polynomial which has $0$ as a repelling fixed point,  and another which has zero as a superattracting fixed point under iteration.  In such a setting,  the filled in Julia set of the semigroup generated by these two polynomials contains a neighbourhood of $0$ and the preimages of $0$ are dense in it.

\medskip

Let $U \subseteq \IC$ be open, and ${\cal P}$ a family of analytic functions, each of which map $U$ into itself. We define the {\em
Julia set} of ${\cal P}$ to be the $z \in U$ such that ${\cal P}$ is not a normal family in any neighbourhood of $z$. 

\begin{lemma}
Let $K \subseteq U \subseteq \IC$ with $K$ compact set and $U$ open. Suppose that 
$\{p_i\;|\; i \in I\}$  is a family of analytic functions which map $U$ into itself, and let ${\cal P}$ be the semigroup that this family generates under composition. Suppose that, for each $z \in K$ there exists a $p \in {\cal P}$ such that
\begin{enumerate}
\item $p(z) \in {\rm int}(K)$
\item $|p'(z)|>L^{-1}$, where $L$ is Landau's constant ($0.5< L \leq 0.5432\ldots$).
\end{enumerate}
Then $\overline{{\rm int}(K)}$ is in the Julia set of ${\cal P}$.
\end{lemma}
Note that Landau's constant is defined as follows:  if $f:\ID\to\IC$ is holomorphic and $f'(0)=1$,  then $f(\ID)$ contains a disk of radius $L$.

\medskip
\noindent{\bf Proof.} 
By 1. and 2., each $z \in K^\circ$ has an open neighbourhood which is mapped into $K^\circ$ by some $p \in {\cal P}$, and on which $|p'(w)|>L^{-1}$. By compactness, a finite set $\{U_1, U_2, \ldots U_n\}$ of such neighborhoods covers $K$. Using compactness of $\overline{U}_i$, we have for each $i$, an $\eta_i>1$, 
such that $p'_i>\eta_i/L$ on $U_i$. Let $r$ be the Lebesgue number of this covering, and $\eta$ the smallest of the $\eta_i$, then, if $z \in K^\circ$ and $D$ is a disc in $K^\circ$ centred at $z$ with radius $s\leq r$,  we
have a $p \in {\cal P}$, such that $p(D) \subseteq K^\circ$ and
contains a disc of radius $\eta s$. By iterating this process we can apply a succession of functions in ${\cal P}$ to $D$, such that was successive images of $D$ contain a disks of radius $s, \eta s, \eta^2 s \ldots$ up to a radius of $r$, which is thereafter maintained. Since $s$ can be chosen arbitrarily small, it follows that ${\cal P}$ cannot be a normal family on any open subset of $K^\circ$. \hfill $\Box$

\begin{lemma}
Let $p$ and $q$ be entire functions, with a common fixed point $c$, which is superattractive for $p$ and repulsive for $q$. Let ${\cal P}$ be the semigroup $\langle p,q \rangle$,
then the Julia set of ${\cal P}$ contains a neighbourhood of $c$. 
\end{lemma}
\noindent{\bf Proof.} This is a standard ``push me,  pull you'' argument which we sketch.  We may assume that $c=0$. We have $p(z)=az^m+O(z^{m+1})$, $q(z)=\mu(z+bz^2+O(z^3))$, where $a\neq 0$, $m \geq 2$ and $|\mu|>1$. We construct a sequence of functions $\{f_n\}$ inductively by $f_0(z)=z$ and
$f_{n+1}(z)=g(f(z))$, where $g$ is either $p$ or $q$. We choose $r>0$ to be sufficiently small that we can ignore higher
degree terms in $p$ and $q$. Let $z_0$ be chosen with $|z_0|<r$, suppose that $f_0,f_1, \ldots f_k$ have already been defined, and set $z_i=f_i(z_0)$. If $|z_k| \geq r$, then let $f_{n+1}(z)=p(f(z))$; otherwise let $f_{n+1}(z)=q(f(z))$.
As soon as $|z_k|\geq r$, the next number $z_{k+1}$ is much smaller; then the $z_i$
gradually increase in size (because $|\mu|>1$), until eventually it exceeds $r$, and the process begins again. The sequence
$\{z_i\}$ is bounded above and below $|a|r^m\leq z_i \leq |\mu|r$.

If $f_{n+1}(z)=p(f_n(z))$, then the logarithmic derivative
\[
\frac{f'_{n+1}(z)}{f_{n+1}(z)}=p'(z_n)\frac{z_n}{p(z_n)} \frac{f'_n(z)}{f_n(z)} \simeq m\frac{f'_n(z)}{f_n(z)}
\]
If $f_{n+1}(z)=q(f_n(z))$, then
\[
\frac{f'_{n+1}(z)}{f_{n+1}(z)}=q'(z_n)\frac{z_n}{q(z_n)} \frac{f'_n(z)}{f_n(z)} \simeq \frac{z+2bz^2}{z+bz^2}\frac{f'_n(z)}{f_n(z)}\simeq (1+bz) \frac{f'_n(z)}{f_n(z)}
\]
If $|z_{n-1}| \geq r$, then $|z_n| \geq |a|r^m$, and it takes $t$ applications of $q$ to get the size of $z_i$ over $r$ again, where $t$ is at most about $\log_{|\mu|} (1/(|a|r^{m-1}))$, in the course of which we multiply the absolute value of the logarithmic derivative by at least
\[
(1-|bz_n|)(1-|bz_{n+1}|) \ldots (1-|bz_{n+t}|) \geq {(1-|b|r)}^t \geq {(1-|b|r)}^{C\log(1/r)}
\]
which can be made as close to 1 as we like  by taking $r$ sufficiently small.
Each time we apply $p$, we multiply the logarithmic derivative by approximately $m$. It follows that $|\frac{f'_n(z)}{f_n(z)}| \to \infty$ as $n \to \infty$. Since the $|f_n(z)|$ is bounded below, it also follows that
$|f'_n(z)|\to \infty$ as well. Thus no subsequence of $\{f_n(z)\}$ can converge to an analytic function. Since
$|f_n(z)|$ is bounded above, $\{f_n(z)\}$ cannot converge to $\infty$ either. Thus $\langle p,q \rangle$ is not a normal
 family on any neighbourhood of $z$. \hfill $\Box$
 
 \bigskip
 
To complete our proof we recall from \S \ref{examples} the  trace polynomials 
$\gamma  (\gamma -\beta )$ from the word $ bab $,  
$ (\beta -\gamma +1)^2 \gamma $  from the word $ babab$,  $\gamma (1-2 \beta +\gamma ^2-(\beta -2) \gamma )$ from the word $ baba^{-1} b$,  $  \gamma(\gamma-\beta) (\beta -\gamma +2)^2 $ from the word  $bababab$\\
 and 
$\gamma ^3(\gamma -\beta ) (\beta +4)  (\beta  (\gamma ^2-3\gamma-4)-\beta ^2( \gamma+1) +4 \gamma ^2+4 \gamma +1)$ from the word $ba^{-2}bababa^{-2}bab$. 

The last polynomial here  is superattractive at $z=0$ and the rest have multipliers at $0$ of $-\beta$, ${(1+\beta)}^2$, $1-2 \beta$ and ${(1-3\beta)}^2$ respectively, so that for each $\beta$ at least one of them has $z=0$ as a repulsive fixed point. Thus, by the lemma, the Julia set of the trace polynomials contains a neighbourhood of $0$, and so a zero-free region $U$ has an image under a trace polynomial into a region which intersects the Julia set, contradicting the fact that these polynomials generate a normal family on $U$.


\begin{thebibliography}{99}  

\bibitem{ASWY} H. Akiyoshi, M. Sakuma, M. Wada and Y. Yamashita, {\em  Punctured torus groups and two bridge knot groups (I)},  Lecture Notes in Mathematics 1909, Springer-Verlag Berlin Heidelberg, 2007.
 
\bibitem{B} A. Beardon, {\em The geometry of discrete groups},
Springer--Verlag, 1983.


\bibitem{Cao} C. Cao,  {\em Some trace inequalities for discrete groups of M\"obius transformations}, Proc. Amer. Math. Soc., {\bf 123}, (1995), 3807--3815.

 

\bibitem{conder} M.D.E. Conder, G.J. Martin and A.
Torstensson, {\em  Maximal symmetry groups of hyperbolic 3-manifolds}, New Zealand J. Math., {\bf 35}, (2006), 37--62.

\bibitem{FK} R. Fricke and F. Klein, {\em Vorlesungen \"uber die Theorie der automorphen Functionen}, Chapter 2, Teubner, Leipzig, 1897.

\bibitem{Gab} D. Gabai, {\em On the Geometric and Topological Rigidity of Hyperbolic 3-Manifolds}, J. American Math. Soc.,
{\bf 10}, (1997),  37--74.

\bibitem{GMT} D. Gabai, R. Meyerhoff and N. Thurston, {\em Homotopy
hyperbolic 3-manifolds are hyperbolic},
Ann. of Math., {\bf 157}, (2003),  335--431.

\bibitem{GMMR} F. W. Gehring, C.Maclachlan, G. J. Martin and A. W. Reid {\em
Arithmeticity, Discreteness and Volume}, Trans. Amer. Math. Soc., {\bf 349},
(1997), 3611 -- 3643.

\bibitem{GMit}  F. W. Gehring and G. J. Martin, {\em Iteration theory and inequalities for Kleinian groups}, Bull. Amer. Math. Soc., {\bf 21}, (1989),   57--63.

\bibitem{GMHalmos} F. W. Gehring and G. J. Martin, {\em Some universal constraints for discrete M\"obius groups}, Paul Halmos; Celebrating 50 Years
of Mathematics, Springer-Verlag, New York, pp. 205--220.

\bibitem{GM3} F. W. Gehring and G. J. Martin, {\em Commutators,
collars and the geometry of M\"obius groups },
J. d'Analyse Math., {\bf 63}, (1994), 175--219.

\bibitem{GMMargulis} F. W. Gehring and G. J. Martin, {\em (p,q,r)-Kleinian groups and the Margulis constant},  Complex analysis
and dynamical systems II, Contemp. Math., 382, (2005), 149--169.

\bibitem{GMannals} F.  W.
Gehring and G. J. Martin, {\em Minimal covolume lattices {\bf I}: spherical points of a  Kleinian group},  Annals of Math., {\bf 170}, (2009), 123--161.

\bibitem{HM1} A. Hinkkanen and G. J. Martin, {\em The Dynamics of Semigroups of Rational Functions I},
Proc. London Math. Soc., {\bf  73}, (1996), 358--384.


\bibitem{HM2} A. Hinkkanen and G. J. Martin, {\em Julia sets of rational semigroups},   Math. Z. {\bf 222},  (1996), 161--169.

\bibitem{Horo} R. Horowitz, {\em Characters of free groups represented in the two-dimensional linear group},  Comm. Pure Appl. Math., {\bf 25}, (1972), 635--649.

\bibitem{Jorg} T. J\o rgensen, {\em On discrete groups of Mobius transformations}, Amer. J. Math., {\bf 98},  (1976),   739--749.

\bibitem{KS}  L. Keen and C. Series, {\em The Riley Slice of Schottky Space},
Proc. London Math. Soc., {\bf 69}, (1994), 72--90.

\bibitem{MM} T. H. Marshall and G.J. Martin, {\em Volumes of
hyperbolic 3-manifolds: Notes on a paper of D. Gabai, G. Meyerhoff
and P. Milley}, J. Conf. Geom. and Dynamics, {\bf 7}, 34--48, 2003.

\bibitem{MM2} T. H. Marshall and G.J. Martin, {\em Minimal co-volume hyperbolic lattices, II: Simple torsion in a Kleinian group}, Ann. Math., {\bf 176}, (2012), 261--301.

\bibitem{MR} C. Maclachlan and A. Reid,  {\em The arithmetic of hyperbolic 3-manifolds,} Springer--Verlag, {\bf 219}, 2003.

\bibitem{SS} R. Stankewitz and H. Sumi, {\em Dynamical properties and structure of Julia sets of postcritically bounded
polynomial semigroups},    Trans. Amer. Math. Soc., {\bf 363}, (2011),  293--5319.

\bibitem{NS} H. Namazi and J. Souto, {\em Non-realizability and ending laminations: Proof of the density conjecture,} Acta Math., {\bf 209},  (2012), 323--395.


\bibitem{sumi1} H. Sumi,  {\em On dynamics of hyperbolic rational semigroups}, J. Math. Kyoto Univ. 37 (1997) 717--733.

\bibitem{sumi2}  H. Sumi,  {\em Rational semigroups, random complex dynamics, and singular functions on the complex plane},  SUGAKU {\bf 61},  (2009), 133--161.

\bibitem{Traina} C. R. Traina, {\em Trace polynomial for two generator subgroups of  $SL(2,\IC)$},  Proc. Amer. Math.Soc.,  {\bf 79}, (1980), 369--372.

\bibitem{V} John Voight, {\em Quaternion Algebras},  unpublished lecture notes, available at
https://math.dartmouth.edu/~jvoight/quat-book.pdf
\end{thebibliography}
\end{document}